\documentclass[preprint,12pt]{elsarticle}

\usepackage{amssymb}
\usepackage{blindtext} % blindtext
\usepackage{booktabs}
\usepackage{eurosym}
\usepackage{siunitx}
\DeclareSIUnit{\EUR}{\mbox{\euro}}
\DeclareSIUnit{\USD}{\mbox{\$}}
\DeclareSIUnit{\year}{\mbox{y}}
\DeclareSIUnit{\betaCo}{\mbox{m$^{2\,\beta}$}}
\usepackage{mathrsfs}
\usepackage{float}
\usepackage{threeparttable}
\usepackage{amsmath}

\usepackage{nomencl}
\usepackage{etoolbox}
\usepackage{svg}
\usepackage{graphicx}
\usepackage{float}
\usepackage{comment}
\usepackage{multirow}
\usepackage[justification=centering]{caption}
\usepackage[super]{nth}

\usepackage{verbatim}

\usepackage{xargs}    % Use more than one optional parameter in a new commands
\usepackage[colorinlistoftodos,prependcaption,textsize=small]{todonotes} %,textsize=tiny
\newcommandx{\david}[2][1=]{\todo[linecolor=red,backgroundcolor=red!25,bordercolor=red,#1]{#2{\\\hfill\tiny{David}}}}
\newcommandx{\felix}[2][1=]{\todo[linecolor=blue,backgroundcolor=blue!25,bordercolor=blue,#1]{#2{\\\hfill\tiny{Felix}}}}
\newcommandx{\rene}[2][1=]{\todo[linecolor=Yellow,backgroundcolor=Yellow!25,bordercolor=Yellow,#1]{#2{\\\hfill\tiny{René}}}}

\journal{Computers \& Chemical Engineering}

\begin{document}

\begin{frontmatter}



\title{HENS unchained: MILP implementation of multi-stage utilities with stream splits, variable temperatures and flow capacities.}

% Busting unnecessary utility restrictions.
% Novel perspective on stream and utility definition
% Optimized Waste Heat Utilization using Linearized Heat Exchanger Network Synthesis
% Optimized Waste Heat Utilization using Linearized Heat Exchanger Network Synthesis


\author[inst1]{David Huber}
\author[inst1]{Felix Birkelbach}
\author[inst1]{René Hofmann}

\affiliation[inst1]{
            organization={Institute for Energy Systems and Thermodynamics, TU Wien},%Department and Organization
            addressline={Getreidemarkt 9/BA}, 
            city={Wien},
            postcode={1060}, 
            country={Austria}}

\begin{abstract}
Heat exchanger network synthesis (HENS) is a well-studied method in research for determining cost-optimal heat exchanger networks. In this paper, we present a modified superstructure formulation to implement streams with variable temperatures and flow capacities. To apply fast MILP solvers, all nonlinear terms, such as those of LMTD, HEX areas and energy balances, are piecewise-linear approximated with simplex or hyperplane models. The translation to MILP is achieved with highly efficient logarithmic coding. One promising application is implementing utilities as streams with variable temperatures and flow capacities. On the one hand, this enables multi-stage heat transfer with stream splits and intermediate utility placement. On the other hand, the temperatures of the utilities can be included as a design parameter in optimizing the heat exchanger network. This makes sense if only the sensible heat of, e.g., thermal oil, water or flue gas, is used as a utility where the inlet and outlet temperatures do not necessarily have to be specified a priori. To examine whether the implementation of utilities as streams leads to more cost-effective solutions, three representative case studies were considered. The results show that reducing the outlet temperature of cold utilities or increasing the outlet temperature of hot utilities leads to significant cost savings. We show that implementing utilities as multi-staged streams with stream splits, variable temperatures and flow capacities is a highly efficient tool for indirect, cost-efficient utility design.
\end{abstract}


\begin{highlights}
\item Multi-stage utilities with stream splits, variable temperatures and flow capacities.
\item Utility temperature targeting enables cost-efficient process integration.
\item Piecewise linear approximation of all non-linear functions: reduced heat exchanger areas, LMTD and energy balances.
\item Linearized problem is solvable with fast MILP solvers.
\end{highlights}


\end{frontmatter}


%% main 
\makenomenclature

% creat sets
\renewcommand\nomgroup[1]{%
  \item[\bfseries
  \ifstrequal{#1}{A}{Acronyms}{%
  \ifstrequal{#1}{P}{Parameters}{%
  \ifstrequal{#1}{V}{Variables}{%
  \ifstrequal{#1}{L}{Subscripts}{%
  \ifstrequal{#1}{H}{Superscripts}{%
  \ifstrequal{#1}{S}{Sets}}}}}}%
]}

% increas spacing
\setlength{\nomlabelwidth}{1.5cm}

% add the units
\newcommand{\nomunit}[1]{%
\renewcommand{\nomentryend}{\hspace*{\fill}#1}}

\nomenclature[P]{$c_{\textrm{cu}}$}{cost coefficient for cold utilities \nomunit{\SI{}{\EUR\per{\kilo\watt\year}}}}

\nomenclature[P]{$c_{\textrm{hu}}$}{cost coefficient for hot utilities \nomunit{\SI{}{\EUR\per{\kilo\watt\year}}}}

\nomenclature[P]{$c_{\textrm{cs}}$}{cost coefficient for cold streams \nomunit{\SI{}{\EUR\per{\kilo\watt\year}}}}

\nomenclature[P]{$c_{\textrm{hs}}$}{cost coefficient for hot streams \nomunit{\SI{}{\EUR\per{\kilo\watt\year}}}}

\nomenclature[L]{$\textrm{cu}$}{cold utility}

\nomenclature[L]{$\textrm{hu}$}{hot utility}

\nomenclature[L]{$\textrm{s}$}{stream}

\nomenclature[L]{$i$}{hot stream}

\nomenclature[A]{US}{utility stream}

\nomenclature[A]{CUS}{cold utility stream}

\nomenclature[A]{HUS}{hot utility stream}

\nomenclature[A]{SSE}{sum of squares error}

\nomenclature[A]{RMSE}{root-mean-square error}

\nomenclature[A]{CS}{case-study}

\nomenclature[A]{GA}{genetic algorithm}

\nomenclature[A]{TAC}{total annual cos}

\nomenclature[L]{$j$}{cold stream}

\nomenclature[L]{$k$}{temperature stage}

\nomenclature[V]{$q$}{heat flow\nomunit{\SI{}{\kilo\watt}}}

\nomenclature[P]{$\Gamma$}{upper bound for temperature difference\nomunit{\SI{}{\celsius}}}

\nomenclature[V]{$F$}{flow capacity\nomunit{\SI{}{\kilo\watt\per\kelvin}}}

\nomenclature[V]{$z$}{binary variable for existance of HEX}

\nomenclature[P]{$c_\textrm{f}$}{step-fixed HEX cost coefficient\nomunit{\SI{}{\EUR\per\year}}}

\nomenclature[P]{$c_\textrm{v}$}{variable HEX cost coefficient\nomunit{\SI{}{\EUR\per{\betaCo\year}}}}

\nomenclature[P]{$\beta$}{cost exponent\nomunit{}}

\nomenclature[P]{$n$}{dimension\nomunit{}}

\nomenclature[P]{$\omega$}{lower bound for heat exchange\nomunit{\SI{}{\kilo\watt}}}

\nomenclature[P]{$\Omega$}{upper bound for heat exchange\nomunit{\SI{}{\kilo\watt}}}

\nomenclature[P]{$U$}{heat transfer coefficient for matches \nomunit{\SI{}{\kilo\watt\per{\meter\squared\kelvin}}}}

\nomenclature[V]{$\mathit{LMTD}$}{logarithmic mean temperature difference
\nomunit{\SI{}{\celsius}}}

\nomenclature[A]{HEN}{heat exchanger network}

\nomenclature[A]{HENS}{heat exchanger network synthesis}

\nomenclature[A]{MILP}{mixed-integer linear programming}

\nomenclature[A]{PtL}{Power-to-Liquid}

\nomenclature[P]{$F$}{flow capacity\nomunit{\SI{}{\kilo\watt\per\kelvin}}}

\nomenclature[P]{$h$}{heat transfer coefficient\nomunit{\SI{}{\kilo\watt\per{\meter\squared\kelvin}}}}

\nomenclature[H]{v}{variable utility paramters}

\nomenclature[P]{$\Delta T_{\textrm{min}}$}{minimum approach temperature
\nomunit{\SI{}{\celsius}}}

\nomenclature[V]{$\Delta T$}{{temperature difference}
\nomunit{\SI{}{\kelvin}}}

\nomenclature[V]{$\mathit{TAC}$}{{total annual costs}
\nomunit{\SI{}{\EUR}}}

\nomenclature[V]{$T$}{{temperature}
\nomunit{\SI{}{\celsius}}}

\nomenclature[H]{in}{inlet}

\nomenclature[H]{out}{outlet}

\nomenclature[A]{HEX}{heat exchanger}

\nomenclature[P]{$N_{\textrm{cs}}$}{number of cold streams}

\nomenclature[P]{$N_{\textrm{hs}}$}{number of hot streams}

\nomenclature[P]{$N_{\textrm{st}}$}{number of stages}

\nomenclature[S]{$HP = \{i \mid i \text{ is hot process stream}; i = 1, \dots N_{\textrm{hs}} \}$}{}

\nomenclature[S]{$CP = \{j \mid j \text{ is cold process stream}; j = 1, \dots N_{\textrm{cs}} \}$}{}

\nomenclature[S]{$ST = \{k \mid k \text{ is temperature stage}; k = 1, \dots N_{\textrm{st}} \}$}{}

\nomenclature[S]{$H_{\textrm{F}} = \{i \mid i \text{ is hot utility stream with variable flow capacity}\}$}{}

\nomenclature[S]{$H_{\textrm{Tin}} = \{i \mid i \text{ is hot utility stream with variable inlet temperature}\}$}{}

\nomenclature[S]{$H_{\textrm{Tout}} = \{i \mid i \text{ is hot utility stream with variable outlet temperature}\}$}{}

\nomenclature[S]{$C_{\textrm{F}} = \{j \mid j \text{ is cold utility stream with variable flow capacity}\}$}{}

\nomenclature[S]{$C_{\textrm{Tin}} = \{j \mid j \text{ is cold utility stream with variable inlet temperature}\}$}{}

\nomenclature[S]{$C_{\textrm{Tout}} = \{j \mid j \text{ is cold utility stream with variable outlet temperature}\}$}{}

\printnomenclature

%\nomenclature[S]{$$}{}
\section{Introduction}
\label{sec:introduction}
Agreed climate targets can only be met through a radical reduction in greenhouse gas emissions. Energy-intensive industries are responsible for a significant share of these emissions. The economic pressure reducing emissions requires companies to pursue cost-effective solutions. One way to reduce emissions cost-effectively is to reduce the energy required to heat or cool process streams. A heat exchanger network (HEN) enables heat exchange between hot and cold process streams. A cost-optimized HEN can be calculated by applying heat exchanger network synthesis (HENS). HENS has a significant impact on energy demand and total annual cost (TAC).

The HEN design problem was first mentioned by Ten Broeck in 1944 \cite{broeck_economic_1944}. The first formal definition was published by Masso \& Rudd in 1969 \cite{masso_synthesis_1969}. All these approaches are sequential methods that decompose the HENS problem into a set of subproblems. Decomposition requires parameter estimation and iterative optimization, which is why global optimality is challenging to achieve. Fully simultaneous methods calculate the optimal utility consumption, stream matches and HEN configuration simultaneously \cite{ciric_heat_1991}. The first simultaneous HENS were published by Yuan et al. \cite{yuan_experiments_1989} in 1989, Yee \& Grossmann \cite{yee_simultaneous_1990} in 1990 and Ciric \& Floudas \cite{ciric_heat_1991} in 1991. For a more detailed elaboration of the historical development, the papers by Furman \& Nikolaos \cite{furman_critical_2002} and Escobar \& Trierweiler \cite{escobar_optimal_2013} are recommended. The latter have shown in their work that Yee \& Grossmann's stage-wise superstructure formulation \cite{yee_simultaneous_1990} gives better results in terms of TAC at lower computation times. Therefore, in this paper, we will build on this formulation.

In Yee \& Grossmann's formulation, two assumptions are made that may lead to sub-optimal HENs: First, implementing utilities is only possible with predefined inlet and outlet temperatures. This assumption is only reasonable where utilities condense a medium at a constant temperature and pressure. If only the sensitive heat in, for example, flue gas or cooling water is used, the temperatures to which the medium must be cooled or heated are of minor importance. Usually, there is a margin for utility temperatures in terms of regulatory and process requirements. Secondly, the utilities must always reach the set temperature in only one heat exchanger without stream splits. In contrast, hot and cold process streams can reach their set target temperature using multi-staged heat exchangers with stream splits. These two limitations inhibit the field application. Considering multi-stage utilities with variable temperatures is essential to optimally integrate the heat sink and source into the process. To run HENS without these assumptions, the Yee \& Grossmann formulation has to be adapted.

Yee \& Grossmann's non-linear formulation belongs to the class of $\mathscr{NP}$-hard problems \cite{furman_computational_2001}. Even with state-of-the-art computational power and solvers, the optimal heat integration of complex industrial processes cannot be calculated. Implementing utilities as streams with variable temperatures and flow capacities further increases the complexity of the optimization problem. Moreover, the non-linear formulation can never guarantee global optimality. Piecewise-linear approximation of the non-linear terms (mean logarithmic temperature difference (LMTD), heat exchanger areas and energy balances) is necessary to find a global minimum within feasible computation time, even though the problem is still $\mathscr{NP}$-hard. Beck \& Hofmann \cite{beck_novel_2018} linearized the superstructure formulation and applied mixed-integer linear programming (MILP) to solve the problem. Compared to the non-linear model, they achieved better results in terms of TAC with shorter computation times.

\subsection{Paper Organization}
\label{sec:introduction_organization}
This paper presents a novel piecewise-linear implementation of utilities as multi-stage streams with stream splits, variable temperatures and flow capacities. The methods in Section \ref{sec:methods} are divided into two main sections. First, all essential adaptations of the superstructure formulation by Yee \& Grossmann are presented in Section \ref{sec:ModificationSuperstructure}. In Section \ref{sec:PiecewiseLinearApprox}, the piecewise-linear approximation of the non-linear terms with hyperplanes and simplices and the transfer to MILP is shown. Section \ref{sec:usecase} introduces three representative use cases from the literature and industrial problems. For each use case, either a cold utility or a cold and a hot utility is implemented as a stream with variable outlet temperature and heat capacity flow. A comparison is made for the results with and without variable utility definitions. We show that minor variations in the utility outlet temperature lead to a significant improvement in terms of TAC. We therefore conclude in Section \ref{sec:conclusion} that variable outlet temperatures and flow capacities allow the cost-optimal design of the necessary utilities.
\section{Methods}
\label{sec:methods}

\subsection{Modification of the Superstructure}
\label{sec:ModificationSuperstructure}
One way of realizing multi-stage utilities is to implement them as streams with stream splits. However, one consequence is that the flow capacity must be specified. If one degree of freedom is blocked by setting the flow capacity for the utility stream, the utilities may not necessarily provide the required energy for heating or cooling the streams. Introducing an additional variable for the flow capacity makes the optimization problem non-linear again. Implementing a variable outlet temperature requires another variable. Both variables are independent and form non-linear relationships, further increasing the complexity of the problem. 

Referring to the stage-wise superstructure according to Yee \& Grossmann \cite{yee_simultaneous_1990}, cold (UC) and hot utilities (UH) can only be located at the end of the streams. The streams can exchange heat in $N_{\textrm{st}}$ stages. 

This paper extends the formulation to implement hot and cold utilities as streams - hereafter referred to as utility streams (US). The objective function 
\begin{equation}
\begin{split}
    \min\mathit{TAC} = \\
    \underbrace{ \sum_{i}{c_{\mathrm{cu}} \, q_{\mathrm{cu},i}} + \sum_{i}\sum_{j}\sum_{k}{c_{\mathrm{cs},j} \, q_{ijk}}}_\text{cold utility costs} +
    \underbrace{\sum_{j}{c_{\mathrm{hu}} \, q_{\mathrm{hu},j}} + \sum_{i}\sum_{j}\sum_{k}{c_{\mathrm{hs},i} \, q_{ijk}}}_\text{hot utility costs} \\
    + \underbrace{\sum_{i}\sum_{j}\sum_{k}{c_{\mathrm{f}}\,z_{ijk}} + \sum_{i}{c_{\mathrm{f}} \, z_{\mathrm{cu},i}} + \sum_{j}{c_{\mathrm{f}} \, z_{\mathrm{hu},j}}}_\text{step-fixed investment costs} \\ 
    + \underbrace{\sum_{i}\sum_{j}\sum_{k}{c_\mathrm{v}\left(\frac{q_{ijk}}{U_{ij}\, \mathit{LMTD}_{ijk}} \right)^{\beta}}}_\text{variable HEX stream costs} \\
    + \underbrace{\sum_{i}{c_\mathrm{v} \left(\frac{q_{\mathrm{cu},i}}{U_{\mathrm{cu},i}\, \mathit{LMTD}_{\mathrm{cu},i}} \right)^{\beta}}}_\text{variable HEX cold utility costs}  + \underbrace{\sum_{j}{c_\mathrm{v} \left(\frac{q_{\mathrm{hu},j}}{U_{\mathrm{hu},j}\, \mathit{LMTD}_{\mathrm{hu},j}} \right)^{\beta}}}_\text{variable HEX hot utility costs}
\end{split} 
\label{eq:TAC}
\end{equation}
where
\begin{equation}
\begin{split}
    U_{ij} = \left( \frac{1}{h_{i}} + \frac{1}{h_{j}} \right)^{-1} \quad &: \quad i \in \mathit{HP}, \, j \in \mathit{CP}\\
    U_{\mathrm{cu},i} = \left( \frac{1}{h_{\textrm{cu}}} + \frac{1}{h_{i}} \right)^{-1} \quad &: \quad i \in \mathit{HP}\\
    U_{\mathrm{hu},j} = \left( \frac{1}{h_{\textrm{hu}}} + \frac{1}{h_{j}} \right)^{-1} \quad &: \quad j \in \mathit{CP}
\end{split} 
\end{equation} and
\begin{equation}
\begin{split}
    \mathit{LMTD}_{ijk} = \frac{\Delta T_{i,j,k} - \Delta T_{i,j,k+1}}{\ln{\frac{\Delta T_{i,j,k}}{\Delta T_{i,j,k+1}}}} \quad &: \quad i \in \mathit{HP}, \, j \in \mathit{CP}, \, k \in \mathit{ST} \\
    \mathit{LMTD}_{\textrm{cu},i} = \frac{\Delta T_{\textrm{cu1},i} - \Delta T_{\textrm{cu2},i}}{\ln{\frac{\Delta T_{\textrm{cu1},i}}{\Delta T_{\textrm{cu2},i}}}} \quad &: \quad i \in \mathit{HP} \\
    \mathit{LMTD}_{\textrm{hu},j} = \frac{\Delta T_{\textrm{hu1},j} - \Delta T_{\textrm{hu2},j}}{\ln{\frac{\Delta T_{\textrm{hu1},j}}{\Delta T_{\textrm{hu2},j}}}} \quad &: \quad j \in \mathit{CP} \\
\end{split} 
\label{eq:LMTD}
\end{equation}
minimizes the TAC of the heat exchanger network. 

Implementing US requires the allocation of costs. The objective in Equation \eqref{eq:TAC} is modified so that costs can be assigned to each hot or cold stream by the cost vector $c_{\mathrm{cs}}$ and $c_{\mathrm{hs}}$. 
\begin{equation}
    c_{\mathrm{hs},i} = \begin{cases} 
        c_{\mathrm{hu}} & \text{if stream }i\text{ is hot utility stream} \\
        0 & \text{otherwise}
    \end{cases}
    \quad : \quad i \in \mathit{HP}
\label{eq:chs}
\end{equation}
\begin{equation}
    c_{\mathrm{cs},j} = 
    \begin{cases} 
        c_{\mathrm{cu}} & \text{if stream }j\text{ is cold utility stream} \\
        0 & \text{otherwise}
    \end{cases}
    \quad : \quad j \in \mathit{CP}
\label{eq:ccs}
\end{equation}
The cost vector given by Equation \eqref{eq:chs} and \eqref{eq:ccs} maps the cost of the hot and cold utilities to the utility streams. This formulation implies that the utility costs are proportional to the heat flow and do not depend on the temperature.

If the flow capacity and the outlet temperature are constant, Equation \eqref{eq:constUtilityLoad} is used to constrain the utility heat loads.
\begin{equation}
\begin{split}
    q_{\textrm{cu,}i} = F_i \, \left( T_{i,k=N_{\textrm{st}}+1} - T_{i,k=N_{\textrm{st}}+2}\right) \quad &: \quad i \in \mathit{HP} \textrm{\textbackslash} \left( H_{\textrm{F}} \cap H_{\textrm{Tout}} \right) \\
    q_{\textrm{hu,}j} = F_j \, \left( T_{j,k=1} - T_{j,k=2} \right) \quad &: \quad j \in \mathit{CP} \textrm{\textbackslash} \left( C_{\textrm{F}} \cap C_{\textrm{Tout}} \right) \\
\end{split}
\label{eq:constUtilityLoad}
\end{equation}
If hot utility streams (HUS) and/or cold utility streams (CUS) are implemented, the utilities are no longer necessary and disabled with Equation \eqref{eq:blockUtilities}.
\begin{equation}
\begin{split}
    z_{\textrm{cu,}i}=0 \quad &: \quad i \in H_{\textrm{F}} \\
    z_{\textrm{hu,}j}=0 \quad &: \quad j \in C_{\textrm{F}} 
\end{split}
\label{eq:blockUtilities}
\end{equation}
The heat exchange between utilities and streams always occurs at the stream ends in only one stage and without stream splits. This results in a total of $N_{\textrm{st}}+1$ stages for heat exchange with other streams and the utility. Due to the disabled utilities, only $N_{\textrm{st}}$ stages are available for the US heat exchange. Increasing the number of stages by one ensures that the same number of stages are available for heat exchange compared to the original superstructure formulation.
\begin{equation}
\begin{split}
    z_{i,j,k=1}=0 \quad &: \quad i \in \mathit{HP}  \textrm{\textbackslash} H_{\textrm{F}}, \, j \in \mathit{CP}  \\
    z_{i,j,k=N_{\textrm{st}}+1}=0 \quad &: \quad j \in \mathit{CP} \textrm{\textbackslash} C_{\textrm{F}}, \, i \in \mathit{HP}  \\
\end{split}
\label{eq:blockstreamHE}
\end{equation}
Blocking the stream heat exchange with Equation \eqref{eq:blockstreamHE} at the added stage secures the stream-to-stream heat exchange at the initial $N_{\textrm{st}}+1$ stages.

%Stream-wise energy balance
The temperatures at position \mbox{$k=1$} and \mbox{$k=N_{\textrm{st}}+2$} in Equation \eqref{eq:ebsh} and \eqref{eq:ebsc} correspond to the inlet and outlet temperatures of the streams. 
\begin{equation}
    \sum_{j}\sum_{k} \, q_{ijk} + q_{\textrm{cu,}i} = F_i \, \left( T_{i,k=1} - T_{i,k=N_{\textrm{st}}+2} \right) \quad : \quad i \in \mathit{HP}
\label{eq:ebsh}
\end{equation}
\begin{equation}
    \sum_{i}\sum_{k} \, q_{ijk} + q_{\textrm{hu,}j} = F_j \, \left( T_{j,k=1} - T_{j,k=N_{\textrm{st}}+2} \right) \quad : \quad j \in \mathit{CP}
\label{eq:ebsc}
\end{equation}
If at least two of the three variables on the right side are assigned a discrete value with Equation \eqref{eq:flowCapConst}, \eqref{eq:inletTempConst} or \eqref{eq:outletTempConst}, the constraints of the stream-wise energy balance remain linear. If fewer values are set, the piecewise-linear approximation presented in Section \ref{sec:linStreamEB} is used.

% stage-wise energy balance
The stage-wise energy balance can be constrained with Equation \eqref{eq:ebsth}. If the flow capacity $F_i$ is not set to a predefined value with Equation \eqref{eq:flowCapConst}, the piecewise-linear approximation from Section \ref{sec:linStreamEB} is used.
\begin{equation}
    \sum_{j} \, q_{ijk} = F_i \, \left( T_{i,k} - T_{i,k+1}\right) \quad : \quad  i \in \mathit{HP}, \,  k \in \mathit{ST}
\label{eq:ebsth}
\end{equation}
\begin{equation}
    \sum_{i} \, q_{ijk} = F_j \, \left( T_{j,k} - T_{j,k+1}\right) \quad : \quad j \in \mathit{CP}, \, k \in \mathit{ST}
\label{eq:ebstc}
\end{equation}

%Assignment of flow capacities:
The flow capacities are set to a specific value with Equation \eqref{eq:flowCapConst}. Otherwise, $F$ is bounded to the predefined range $[F_{\textrm{min}}^{\textrm{set}},  F_{\textrm{max}}^{\textrm{set}}]$ with Equation \eqref{eq:flowCapVar}.
\begin{equation}
\begin{split}
    F_{i} = F_{i}^{\textrm{set}} \quad : \quad i \in \mathit{HP} \textrm{\textbackslash} H_{\textrm{F}} \\ F_{j} = F_{j}^{\textrm{set}} \quad : \quad j \in \mathit{CP} \textrm{\textbackslash} C_{\textrm{F}}
\end{split}
\label{eq:flowCapConst}
\end{equation}
%Variable flow capacities can be considered according to Equation \eqref{eq:flowCapVar} using predefined boundaries.
\begin{equation}
\begin{split}
    F_{i,\textrm{min}}^{\textrm{set}} \le F_{i} \le F_{i,\textrm{max}}^{\textrm{set}} & \quad : \quad i \in H_{\textrm{F}} \\
    F_{j,\textrm{min}}^{\textrm{set}}\le F_{j} \le F_{j,\textrm{max}}^{\textrm{set}} & \quad : \quad j \in C_{\textrm{F}}
\end{split}
\label{eq:flowCapVar}
\end{equation}

%Assignment of inlet temperatures:
Constant inlet or outlet temperatures are set with Equations \eqref{eq:inletTempConst} and \eqref{eq:outletTempConst}. Variable temperatures are constrained to a specified range for the inlet temperature $[T^{\textrm{in}}_{\textrm{min}}, T^{\textrm{in}}_{\textrm{max}}]$ and the range for the outlet temperature $[T^{\textrm{out}}_{\textrm{min}}, T^{\textrm{out}}_{\textrm{max}}]$ using Equations \eqref{eq:inletTempVar} and \eqref{eq:outletTempVar}, respectively.
\begin{equation}
\begin{split}
    T_{i,k=1} = T_{i}^{\textrm{in}} \quad &: \quad i \in \mathit{HP} \textrm{\textbackslash} H_{\textrm{Tin}} \\
    T_{j,k=N_{\textrm{st}}+2} = T_{j}^{\textrm{in}} \quad &: \quad j \in \mathit{CP} \textrm{\textbackslash} C_{\textrm{Tin}}
\end{split}
\label{eq:inletTempConst}
\end{equation}

\begin{equation}
\begin{split}
    T_{i,k=N_{\textrm{st}}+2} = T_{i}^{\textrm{out}} \quad : \quad i \in \mathit{HP} \textrm{\textbackslash} H_{\textrm{Tout}} \\
    T_{j,k=1} = T_{j}^{\textrm{out}} \quad : \quad j \in \mathit{CP} \textrm{\textbackslash} C_{\textrm{Tout}}
\end{split}
\label{eq:outletTempConst}
\end{equation}

\begin{equation}
\begin{split}
    T_{i,\textrm{min}}^{\textrm{in}} \le T_{i,k=1} \le T_{i,\textrm{max}}^{\textrm{in}} & \quad : \quad i \in H_{\textrm{Tin}} \\
    T_{j,\textrm{min}}^{\textrm{in}} \le T_{j,k=N_{\textrm{st}}+2} \le T_{j,\textrm{max}}^{\textrm{in}} & \quad : \quad j \in C_{\textrm{Tin}}
\end{split}
\label{eq:inletTempVar}
\end{equation}

\begin{equation}
\begin{split}
    T_{i,\textrm{min}}^{\textrm{out}} \le T_{i,k=N_{\textrm{st}}+2} \le T_{i,\textrm{max}}^{\textrm{out}} & \quad : \quad i \in H_{\textrm{Tout}} \\
    T_{j,\textrm{min}}^{\textrm{out}} \le T_{j,k=1} \le T_{j,\textrm{max}}^{\textrm{out}} & \quad : \quad j \in C_{\textrm{Tout}}
\end{split}
\label{eq:outletTempVar}
\end{equation}

Note that, if stream inlet and outlet temperatures are defined in a specific range, the conditions
\begin{equation}
    \begin{split}
        T_{i,\textrm{min}}^{\textrm{in}} \ge T_{i,\textrm{max}}^{\textrm{out}} & \quad : \quad i \in H_{\textrm{Tin}} \cap H_{\textrm{Tout}}\\
        T_{j,\textrm{max}}^{\textrm{in}} \le T_{j,\textrm{min}}^{\textrm{out}} & \quad : \quad j \in C_{\textrm{Tin}} \cap C_{\textrm{Tout}}
    \end{split}
    \label{eq:tempCondition}
\end{equation}
must always be fulfilled to obtain a feasible solution.

The following constraints are not affected by variable temperatures or flow capacities.

Monotonic decrease in temperature:
\begin{equation}
\begin{split}
    T_{i,k} &\ge T_{i,k+1} \quad : \quad  i \in \mathit{HP}, \,  k \in \mathit{ST} \\
    T_{j,k} &\ge T_{j,k+1} \quad : \quad  j \in \mathit{CP}, \,  k \in \mathit{ST} \\
\end{split}
\label{eq:monotonicDecrease}
\end{equation}

Bounds for heat loads:
\begin{equation}
\begin{split}
    z_{i,j,k} \, \omega_{\textrm{s}} \le q_{i,j,k} \le z_{i,j,k} \, \Omega_{\textrm{s}} \quad &: \quad i \in \mathit{HP}, \, j \in \mathit{CP}, \, k \in \mathit{ST} \\
    z_{\textrm{cu,}i} \, \omega_{\textrm{cu}} \le q_{\textrm{cu,}i} \le z_{\textrm{cu,}i} \, \Omega_{\textrm{cu}} \quad &: \quad i \in \mathit{HP} \\
    z_{\textrm{hu,}j} \, \omega_{\textrm{hu}} \le q_{\textrm{hu,}j} \le z_{\textrm{hu,}j} \, \Omega_{\textrm{hu}} \quad &: \quad j \in \mathit{CP}
\end{split}
\end{equation}

Bounds for temperature differences:
\begin{equation}
\begin{split}
    \Delta T_{\textrm{min}} \le \Delta T_{i,j,k} \\
    \Delta T_{i,j,k} \le T_{i,k} - T_{j,k} + \Gamma \left( 1 - z_{i,j,k} \right) \\
    \Delta T_{i,j,k+1} \le T_{i,k+1} - T_{j,k+1} + \Gamma \left( 1 - z_{i,j,k} \right)
\end{split}
\, : \, i \in \mathit{HP}, \, j \in \mathit{CP}, \, k \in \mathit{ST}
\end{equation}
\begin{equation}
\begin{split}
    \Delta T_{\textrm{min}} \le \Delta T_{\mathrm{cu1},i} \le T_{i,k=N_{\textrm{st}}+1} - T_{\textrm{cu}}^{\textrm{out}} + \Gamma \left( 1-z_{\textrm{cu,}i} \right)\\
    \Delta T_{\textrm{min}} \le \Delta T_{\mathrm{cu2},i} \le T_{i,k=N_{\textrm{st}}+2} - T_{\textrm{cu}}^{\textrm{in}} + \Gamma \left( 1-z_{\textrm{cu,}i} \right)
\end{split}
\quad : \quad i \in \mathit{HP}
\end{equation}
\begin{equation}
\begin{split}
    \Delta T_{\textrm{min}} \le \Delta T_{\mathrm{hu1},j} \le T_{\textrm{hu}}^{\textrm{out}} - T_{j,k=2} + \Gamma \left( 1-z_{\textrm{hu,}j} \right)\\
    \Delta T_{\textrm{min}} \le \Delta T_{\mathrm{hu2},j} \le T_{\textrm{hu}}^{\textrm{in}} -  T_{j,k=1} + \Gamma \left( 1-z_{\textrm{cu,}j} \right)
\end{split}
\quad : \quad j \in \mathit{CP}
\end{equation}

Integrality: 
\begin{equation}
\begin{split}
    z_{i,j,k}, \, z_{\textrm{cu,}i}, \, z_{\textrm{hu,}j} \quad &: \quad i,j,k \in \{0,1\}
\end{split}
\end{equation}

Non-negativity constraints: 
\begin{equation}
\begin{split}
    q_{i,j,k}, \, q_{\textrm{cu,}i}, \, q_{\textrm{hu,}j} \ge 0 \quad &: \quad i \in \mathit{HP}, \, j \in \mathit{CP}, \, k \in \mathit{ST}
\end{split}
\end{equation}

\subsection{Piecewise-linear Approximation}
\label{sec:PiecewiseLinearApprox}
To integrate the design of the utilities into the HENS and find a global optimum within a feasible computation time, piecewise-linear approximation is essential.
For the sake of simplicity, a function with one primary curvature is called convex or concave accordingly. In contrast to the convex heat exchanger area of the streams and the concave heat exchanger surface of the utilities, the energy balance is neither convex nor concave. The energy balance is essentially a multiplication of two independent variables. The resulting saddle-shaped function can no longer be represented with sufficient accuracy by simple concave or convex approximations. Therefore, the following two methods for linear approximation are distinguished in this paper: Piecewise-linear approximation with hyperplanes and with simplices.

Piecewise-linear approximation with hyperplanes is used for the concave function of the stream HEX area (see  Section \ref{sec:streamHEX}) and the concave function of the utility HEX area (see Section \ref{sec:utHEX}). Each hyperplane is defined by a linear function with coefficients $\boldsymbol{a}$, which specifies offset and slope. The coefficients $\boldsymbol{a}$ are determined using a nonlinear optimization that minimizes the sum of squares error (SSE) between the linearized planes and the data points. Concave functions can be linearized in the same way by considering the identity \mbox{$\min(y) = -\max (-y)$}. The accuracy of the approximation can be adjusted by adding hyperplanes until a defined root-mean-square error (RMSE) is reached. In contrast to piecewise-linear approximation with simplices, only limited accuracy can be achieved for non-convex or non-concave approximations.

Convex or concave approximation with hyperplanes is only suitable for functions that curve in only one direction. In the natural sciences, however, problems often occur which require a multiplication of optimization variables. For example, the two-dimensional function \mbox{$f(y) = x_1\,x_2$} is saddle-shaped and cannot be approximated convexly or concavely with sufficient accuracy. By contrast, any continuous function can be approximated piecewise-linearly with simplices. In the two-dimensional set, for a grid with $w$ elements, the function \mbox{$f:[0,w]^2 \to \mathbb{R}$} can be divided into triangles \cite{vielma_modeling_2011}. The function $f$ can thus be approximated with piecewise functions linearly within the triangles. In this paper, the $\boldsymbol{J}_1$ union jack triangulation is used. This method requires a grid with the nodes of the triangles in its intersection. A non-linear optimization problem determines the grid points and the plane equations of the triangles by minimizing the SSE. Piecewise-linear approximation with simplices is used for the stream- and stage-wise energy balances (see Section \ref{sec:linStreamEB}) and the LMTD (see Section \ref{sec:LMTD}).

\subsubsection{Stream Heat Exchanger Area}
\label{sec:streamHEX}
\begin{figure}[H]
    \centering 
    \includegraphics[width=1\textwidth]{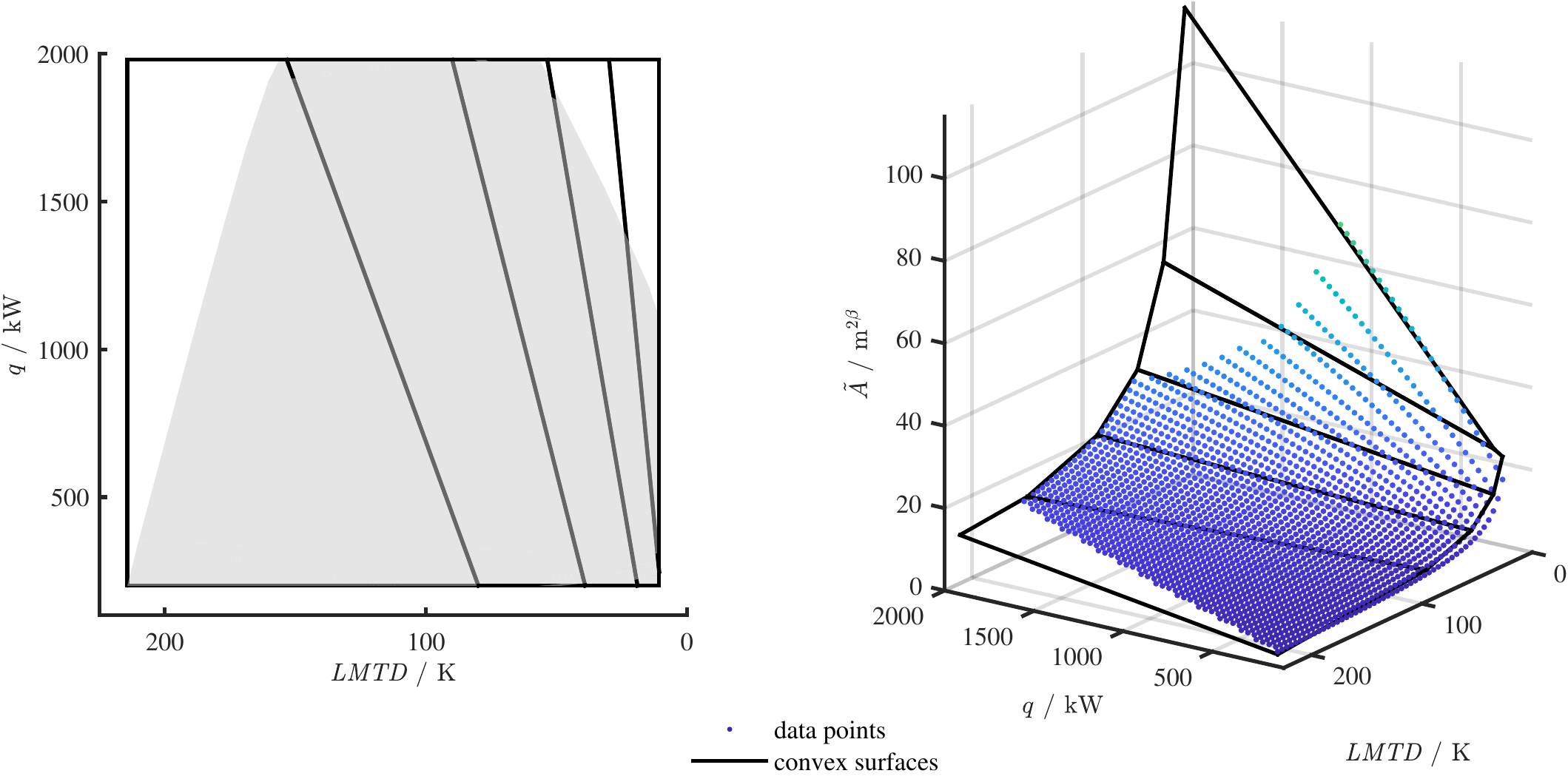}
    \caption{Piecewise-linear approximation of the reduced stream HEX area $\Tilde{A}$ as a function of the heat flow $q$ and $LMTD$ with five hyperplanes. Hot stream: \mbox{$T^{\textrm{in}} = \SI{270}{\celsius}$}, \mbox{$T^{\textrm{out}} = \SI{160}{\celsius}$}, \mbox{}$F = \SI{18}{\kilo\watt\per\kelvin}$, \mbox{$h = \SI{1}{\kilo\watt\per\meter\squared\per\kelvin}$}. Cold stream: \mbox{$T^{\textrm{in}} = \SI{50}{\celsius}$}, \mbox{$T^{\textrm{out}} = \SI{210}{\celsius}$}, \mbox{$F = \SI{20}{\kilo\watt\per\kelvin}$}, \mbox{$h = \SI{1}{\kilo\watt\per\meter\squared\per\kelvin}$}. \mbox{$\beta = 0.8$}. \mbox{$RMSE = \SI{1.26}{\percent}$}.}
    \label{fig:linHEXstream}
\end{figure}

The reduced heat exchanger area 
\begin{equation}
\label{eq:A_HEX}
    \tilde{A}_{ijk} = \left( \frac{q_{ijk}}{U_{ij}\,\textit{LMTD}_{ijk}} \right)^{\beta}
\end{equation}
for a stream HEX is convex. The solution space is reduced to a physically feasible domain as per to Beck et al. \cite{beck_novel_2018}. Beck et al. formulate a linear optimization problem to constrain the independent variables $q$ and $\textit{LMTD}$ to a physically feasible domain. The $N_{\textrm{hyp}}$ hyperplanes of the two-dimensional function are defined with coefficients $\boldsymbol{a}$ such that $\tilde{A} = a_0 + a_1 \, \textrm{LMTD} + a_2 \, q$ for each data point. Figure \ref{fig:linHEXstream} shows the reduced solution space with $2014$ data points in light gray and hyperplane approximation for two example streams. Within this example we are able to achieve an RMSE of $\SI{1.26}{\percent}$ using $5$ hyperplanes. Above $22$ hyperplanes, the RMSE of $\SI{1.16}{\percent}$ does not change within the lsqnonlin solver's step size tolerance of $10^{-6}$.

\subsubsection{Utility Heat Exchanger Area}
\label{sec:utHEX}
\begin{figure}[H]
    \centering
    \includegraphics[width=0.6\textwidth]{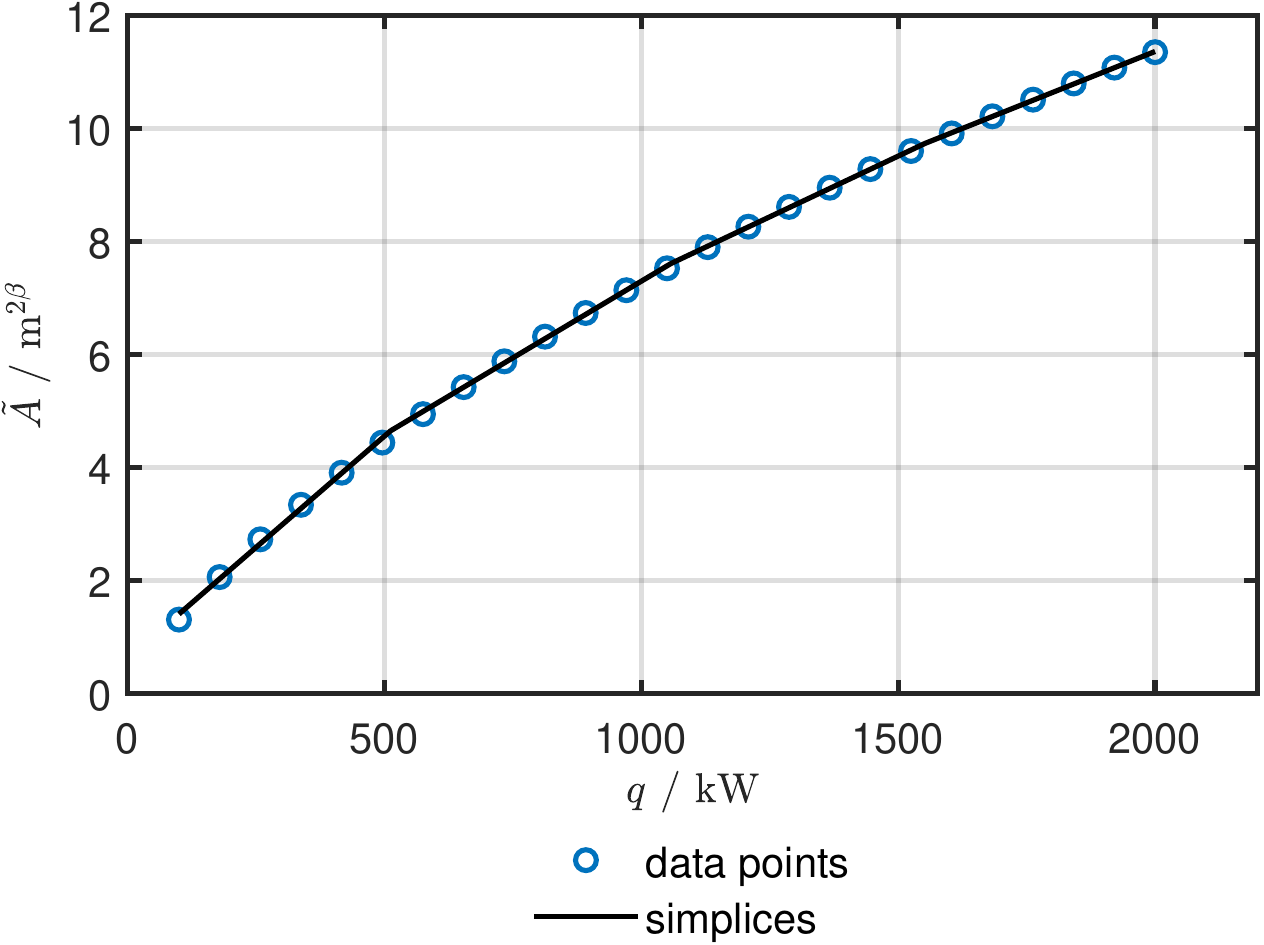}
    \caption{Piecewise-linear approximation of the reduced utility HEX area $\Tilde{A}$ as a function of the heat flow $q$ with four lines. Hot stream: \mbox{$T^{\textrm{in}} = \SI{270}{\celsius}$}, \mbox{$T^{\textrm{out}} = \SI{160}{\celsius}$}, \mbox{$F = \SI{18}{\kilo\watt\per\kelvin}$}, \mbox{$h = \SI{1}{\kilo\watt\per\meter\squared\per\kelvin}$}. Cold utility: \mbox{$T^{\textrm{in}} = \SI{10}{\celsius}$}, \mbox{$T^{\textrm{out}} = \SI{30}{\celsius}$}, \mbox{$h = \SI{1}{\kilo\watt\per\meter\squared\per\kelvin}$}. \mbox{$\beta = 0.8$}. \mbox{$RMSE = \SI{0.38}{\percent}$}.}
    \label{fig:linHEXutility}
\end{figure}

Since three out of four temperatures are fixed at the utility heat exchangers, the reduced heat exchanger area can be formulated as a function of the heat flow $q$ \cite{beck_novel_2018}. The one-dimensional correlation of the reduced heat exchanger area for hot utilities
\begin{equation}
    \label{eq:AhexUH}
    \tilde{A}_{\textrm{hu,}j} \left( q_{\textrm{hu,}j}  \right) = 
    \left[ 
    \frac{q_{\textrm{hu,}j} \, \left( \ln{\left(T_{\textrm{hu}}^{\textrm{in}} - T_j^{\textrm{out}} \right) } - 
    \ln{\left( T_{\textrm{hu}}^{\textrm{out}} - T_{\textrm{j}}^{\textrm{out}} + \frac{q_{\textrm{hu,}j}}{F_j} \right)}\right)
    }{
    U_{\textrm{hu,}j}\,\left( T_{\textrm{hu}}^{\textrm{in}} - T_{\textrm{hu}}^{\textrm{out}} + \frac{q_{\textrm{hu,}j}}{F_j} \right)
    } 
    \right]^{\beta}
\end{equation}
and cold utilities
\begin{equation}
    \label{eq:AhexUC}
    \tilde{A}_{\textrm{cu,}i} \left( q_{\textrm{cu,}i}  \right) = 
    \left[
    \frac{q_{\textrm{cu,}i} \, \left( 
    \ln{\left( T_i^{\textrm{out}} + \frac{q_{\textrm{cu,}i}}{F_i} - T_{\textrm{cu}}^{\textrm{out}} \right)} - \ln{\left( T_i^{\textrm{out}} - T_{\textrm{cu}}^{\textrm{in}}\right) }\right)
    }{
    U_{\textrm{cu,}i}\,\left( \frac{q_{\textrm{cu,}i}}{F_i} - T_{\textrm{cu}}^{\textrm{out}} + T_{\textrm{cu}}^{\textrm{in}}\right)
    } 
    \right]^{\beta}
\end{equation}
is concave. Again, equation \eqref{eq:AhexUH} and \eqref{eq:AhexUC} are restricted to the physically solvable domain and represented by lines. The lines are thus represented as linear equations \mbox{$\tilde{A} = a_0 + a_1 \, q$} for each plane and heat exchanger. The coefficients are again determined by non-linear minimization of the SSE until an RMSE criterion is met. Figure \ref{fig:linHEXutility} shows the concave function of the reduced utility HEX area with 25 data points and the piecewise-linear approximation. With four lines, an RMSE of $\SI{0.38}{\percent}$ can be achieved. In this case, an ideal linear approximation would be possible by interpolating the data points. In this case, the improved accuracy is out of proportion to the required binary variables, which unnecessarily increases the complexity and computation time of the optimization problem.

\subsubsection{Energy Balances}
\label{sec:linStreamEB}
\begin{figure}[H]
    \centering
    \includegraphics[width=1\textwidth]{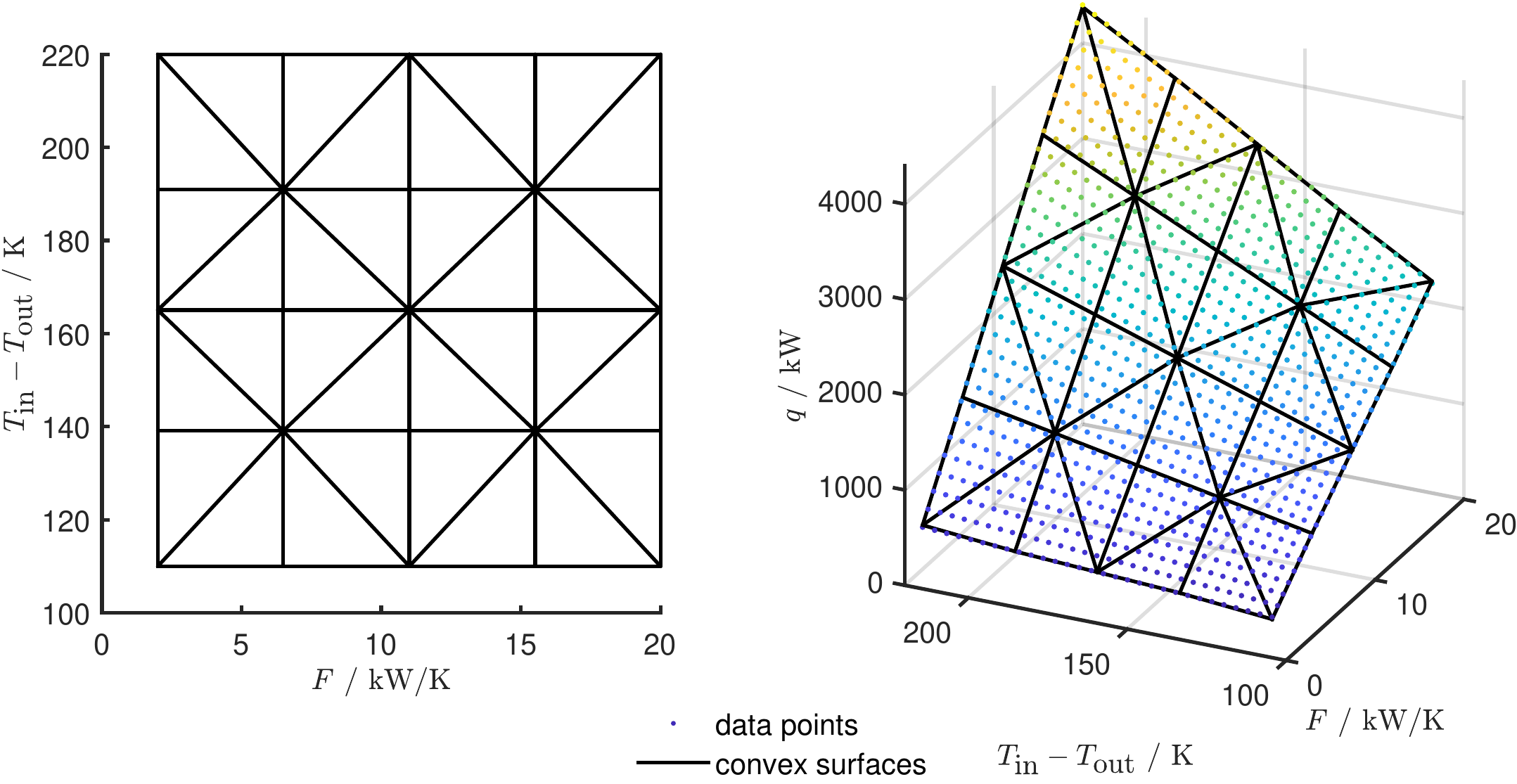}
    \caption{Piecewise-linear approximation of the stream-wise energy balance as a function of the flow capacity $F$ and the temperature difference \mbox{$T_{\textrm{in}}-T_{\textrm{out}}$} with $32$ simplices. Hot stream: \mbox{$T^{\textrm{in}} = \SI{270}{\celsius}$}, \mbox{$T^{\textrm{out}} = [ \SI{50}{},  \SI{160}{}]\,\SI{}{\celsius}$}, \mbox{$F = [ \SI{2}{}, \SI{20}{}]\,\SI{}{\kilo\watt\per\kelvin}$}. \mbox{$RMSE = \SI{0.28}{\percent}$}.}
    \label{fig:energyBalance}
\end{figure}

The piecewise-linear approximation of the energy balances, Equations \mbox{(\ref{eq:ebsh}-\ref{eq:ebstc})}, is of central importance to implementing streams with variable inlet or outlet temperature and flow capacity. Figure \ref{fig:energyBalance} shows the $900$ data points and the piecewise-linear approximation with simplices of a stream-wise-energy balance. The heat flow $q$ is plotted as a function of the flow capacity $F$ and the temperature difference \mbox{$T_{\textrm{in}}-T_{\textrm{out}}$}. The saddle-shaped function is approximated with $32$ simplices on an equidistant 4x4 grid with an RMSE of $\SI{0.28}{\percent}$. 

\subsubsection{LMTD}
\label{sec:LMTD}
\begin{figure}[H]
    \centering
    \includegraphics[width=1\textwidth]{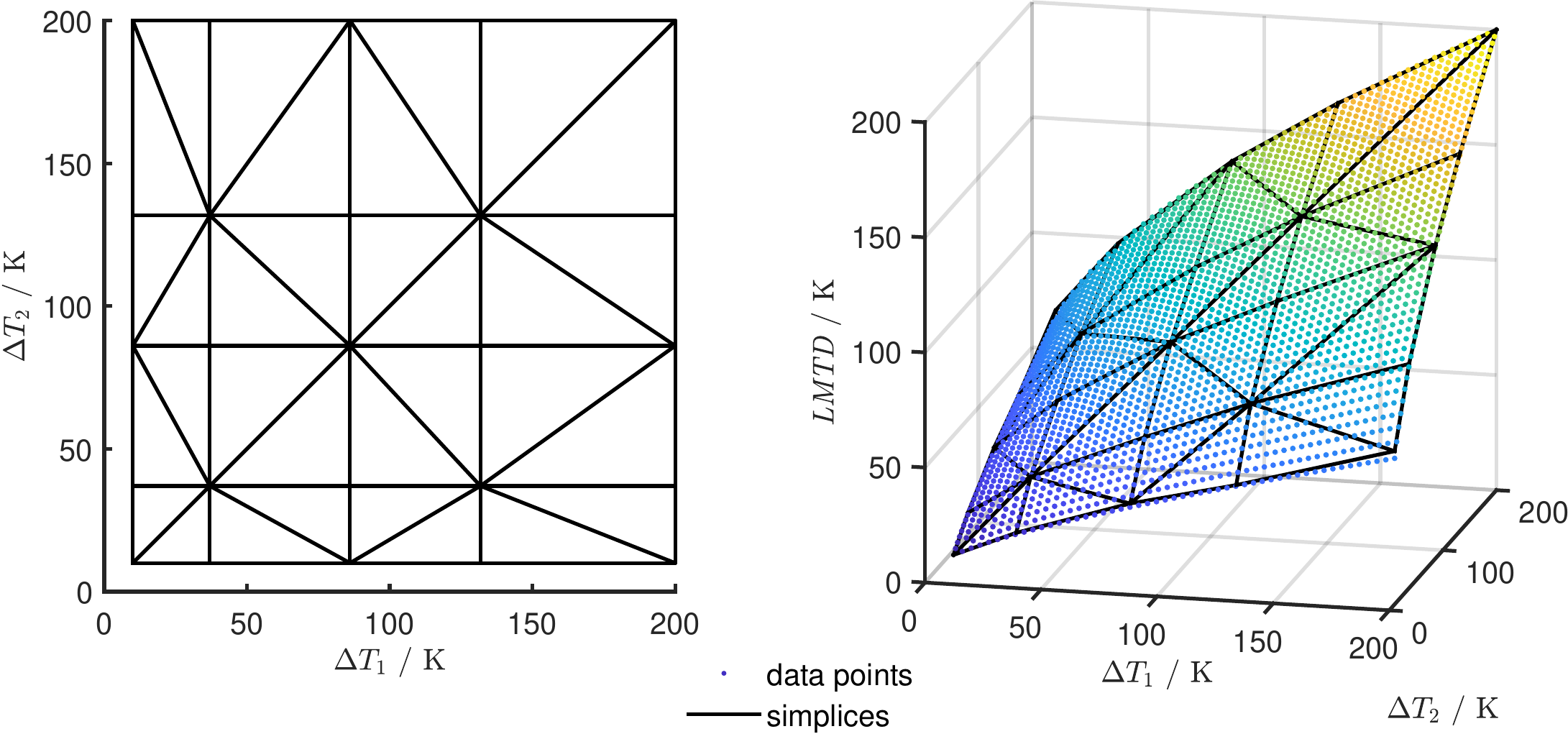}
    \caption{Piecewise-linear approximation of the $LMTD$ as a function of the two temperature differences $\Delta T_1$ and $\Delta T_2$. $LMTD$ in the range from \SI{10}{\kelvin} to \SI{200}{\kelvin}. \mbox{$RMSE = \SI{0.34}{\percent}$}.}
    \label{fig:LMTD}
\end{figure}

The LMTD, according to equation \eqref{eq:LMTD}, is concave and can be approximated with hyperplanes and simplices. Both methods require additional binary variables. The approximation with simplices offers considerable advantages in terms of the MILP translation. Significantly higher accuracies can be achieved with the same number of binary variables. Figure \ref{fig:LMTD} shows the piecewise-linear approximated LMTD with $900$ data points on a 4x4 grid. In regions with larger curvature, more simplices are placed.

\subsection{Translation to MILP}
The translation to MILP should be carried out with as few auxiliary binary variables as possible. Thus, the minimization problem can be solved efficiently within a feasible timeframe.

The streams' convex reduced HEX area is translated most easily to MILP. The hyperplanes shown in Figure \ref{fig:linHEXstream} can be translated to MILP with one inequality each and without additional binary variables, see \cite{beck_novel_2018}. 

For all other functions, binary variables are necessary to translate the simplices into MILP. Vielma \& Nemhauser developed a logarithmic approach in \cite{vielma_modeling_2011} to reduce the number of binary variables. A grid with $w$ elements in an $n$-dimensional space, where $w$ is a power of two, is composed of \mbox{$T = w^n\,n!$} Simplices. The $T$ simplices can be translated to MILP highly efficiently with \mbox{$\lceil \log_2 T \rceil$} binary variables and \mbox{$ \lceil 2 \, \log_2 T \rceil$} additional constraints. The piecewise-linear approximations based on simplices presented in the previous sections are all on a grid with \mbox{$w=4$} elements in each dimension. The one-dimensional approximation of the utility HEX area in Figure \ref{fig:linHEXutility} is modeled with four simplices. Thus, two binary variables and four additional constraints are used to translate the correlation to MILP. On the other hand, the widely used SOS2 approach would require four binary variables. The approximation of the two-dimensional correlations for the streams HEX area, energy balance and LMTD is composed of 32 simplices. These can be translated to MILP with five binary variables and ten additional constraints each. Due to the small number of binary variables combined with the high accuracy of the approximation, non-linear correlations can be approximated highly efficiently and modern MILP solvers can calculate a global optimum in feasible computing time. 

Since not all approximations reach the value \mbox{$f(0) = 0$}, the functions are toggled with binary variables and a \mbox{big-M} approach. The hyperplane approximation's maximum value always occurs in the corner of the domain. Accordingly, the \mbox{big-M} value is chosen. By choosing the smallest possible \mbox{big-M} value, the problem remains tight and the stability of the numerical solving algorithms is improved because the feasible region of the LP relaxation is not unnecessarily expanded \cite{camm_cutting_1990}.

\section{Case studies \& Results}
\label{sec:usecase}
This paper examines whether it is beneficial in terms of TAC to implement utilities as multi-stage streams with stream splits and variable outlet temperature and flow capacity. Furthermore, the influence of the utility outlet temperature on the TAC is studied. For this purpose, in each of the three representative case studies (CS), all utilities that only use sensible heat are implemented as streams. Depending on the utilities, the following cases were considered: 
\begin{description}
    \item[base] For each case study, the base case is used to compare the results with literature values and to validate the optimization framework.
    \item[var UC] A cold utility with variable outlet temperature and flow capacity is implemented when only the sensible heat of a medium such as water or thermal oil is used for cooling. 
    \item[var UC \& UH] A hot and cold utility with variable outlet temperature and flow capacity is implemented when only the sensible heat is used for both cooling and heating.
\end{description}

\subsection{Piecwise-linear Approximation \& Implementation}
Planes were added to the linear models of the stream HEX area until the RMSE was below \mbox{$\SI{1.0}{\percent}$}. To limit the number of binaries used to transfer the simplices to MILP, the approximation of the utility HEX area, energy balances and LMTD were calculated on a \mbox{$4\times4$} grid with 32 Simplices. The $RMSE$ is below \mbox{$\SI{0.5}{\percent}$} for all models.

All optimization problems in this paper have been modeled using Yalmip R20210331 \cite{lofberg_toolbox_2004} in Matlab R2022a. All problems were solved using Gurobi 9.5.2 on a 128-core system (AMD EPYC 7702P) with \mbox{$256$ GB} RAM.

Each optimization problem was solved three times. The solution with the lowest computation time is presented. The convergence behavior over time is shown in \ref{sec:appendixConvergence} with its characteristic values, relative gap, upper and lower objective bounds. The relative gap is defined as the gap between the best feasible solution objective and the best bound. The calculations are terminated if the relative gap is smaller or equal than the tolerance of the MIP solver. The default value is \mbox{$\SI{0.01}{\percent}$}. 

\subsection{Case study 1}
The first case study was presented by Ahmad \cite{ahmad_heat_1985} and is composed of two hot and two cold streams. The stream data is listed in \mbox{Table \ref{tab:UseCase1}}. Since the latent heat of steam is used as hot utility, the outlet temperature of the steam cannot be adjusted without changing the steam parameters. Accordingly, no HUS is implemented. Since only the sensible heat of the cooling water is used, the cold utility is implemented as CUS. The CUS stream definition is marked with the superscript v. The parentheses specify the range of permissible values for the outlet temperature and the flow capacity.

\begin{table}[H]
    \centering
    \caption{Stream data for case study 1: Ahmad \cite{ahmad_heat_1985}.}
    \label{tab:UseCase1}
    \begin{threeparttable}
    \begin{tabular}{crrrr}
        \toprule
        Stream & $T^{\mathrm{in}}$ / $\SI{}{\celsius}$ & $T^{\mathrm{out}}$ / $\SI{}{\celsius}$ & $F$ / $\SI{}{\kilo\watt\per\kelvin}$ & $h$ / $\SI{}{\kilo\watt\per\meter\squared\per\kelvin}$\\
        \midrule 
        H1 & $260$ & $160$ & $3.0$ & $0.4$ \\
        H2 & $250$ & $130$ & $1.5$ & $0.4$ \\
        C1 & $120$ & $235$ & $2.0$ & $0.4$ \\
        C2 & $180$ & $240$ & $4.0$ & $0.4$ \\
        UH & $280$ & $279$ & - & $0.4$ \\
        UC & $30$ & $80$ & - & $0.4$ \\
        UC$^\text{v}$ & $30$ & $[31, 80]$ & $(0, 20]$ & $0.4$ \\
        \bottomrule
    \end{tabular}
    \begin{tablenotes}[param, flushleft]
    \footnotesize
        \item HEX costs: \mbox{$c_{\mathrm{f}} = \SI{0}{\USD\per\year}$}, \mbox{$c_{\mathrm{v}} = \SI{300}{\USD\per\betaCo\per\year}$}, \mbox{$\beta = 0.5$}
        \item Utility costs: \mbox{$c_{\mathrm{hu}} =  \SI{110}{\USD\per\kilo\watt\per\year}$}, \mbox{$c_{\mathrm{cu}} =  \SI{12.2}{\USD\per\kilo\watt\per\year}$}
        \item Min. approach temperature: \mbox{$\Delta T_{\mathrm{min}} = \SI{1}{\celsius}$}
    \end{tablenotes}
    \end{threeparttable}
\end{table}

The number of variables, the computation time, and the relative gap for CS1 can be seen in Table \ref{tab:performanceDataCS1}. The number of binary variables increases significantly with one implemented CUS from $120$ to $376$. Therefore, the computation time for the var UC case of \mbox{$\SI{22.79}{\second}$} is nearly ten times that of the base case with \mbox{$\SI{2.30}{\second}$}.

\begin{table}[htp]
\centering
\caption{Problem size, computation time and relative gap for CS1.}
\label{tab:performanceDataCS1}
\begin{tabular}{lrrrrr}
\toprule
\multicolumn{1}{c}{case} & \multicolumn{1}{c}{variables / -} & \multicolumn{1}{c}{binaries / -}  & \multicolumn{1}{c}{time / s} & \multicolumn{1}{c}{rel. gap / \%} \\
\midrule
 base & $321$ & $120$ & $2.30$ & $0.0000$ \\
 var UC & $3610$ & $376$ & $22.79$ & $0.0009$ \\
 \bottomrule
\end{tabular}
\end{table}

To validate the developed framework, the base case is compared with three references from the literature. The results are summarized in Table \ref{tab:results-UseCase1}. The stream plots of the resulting HENs are shown in \ref{sec:appendixCaseStudyPlotsCS1}. The optimal HENs of Ahmad \cite{ahmad_heat_1985}, Nielsen et al. \cite{nielsen_heat_1996} and Khorasany \& Fesanghary \cite{khorasany_novel_2009} were calculated without using stream splits. Khorasany \& Fesanghary \cite{khorasany_novel_2009} used a two-level approach with a harmony search algorithm and sequential quadratic programming to determine the best known literature value for minimum TAC of \mbox{$1.1895 \cdot 10^4 \, \SI{}{\USD\per\year}$}. In contrast to the literature values, we used three stages instead of four and allowed for stream splits. We were able to find a solution for the base case with \mbox{$1.1792 \cdot 10^4 \, \SI{}{\USD\per\year}$} of TAC. In terms of TAC, the calculated solution is \mbox{$\SI{0.87}{\percent}$} cheaper than the best literature value from Khorasany \& Fesanghary \cite{khorasany_novel_2009}. Since both the values of the TAC and the heat loads show only minor deviations, it can be assumed that the framework presented in this paper provides reliable results. 

\begin{table}[htb]
\centering
\caption{Results for CS1: Comparison of costs and heat loads at the utilities and streams.}
\label{tab:results-UseCase1}
\begin{tabular}{lcccc}
\toprule
\multirow{2}{*}{reference} &  & \multicolumn{3}{c}{heat load / $\SI{}{\kilo\watt}$} \\
 & \multicolumn{1}{c}{TAC / $10^4 \, \SI{}{\USD\per\year}$} & \multicolumn{1}{c}{CU} & \multicolumn{1}{c}{HU} & \multicolumn{1}{c}{stream} \\ \midrule
 Ahmad \cite{ahmad_heat_1985} & $1.2870$ & $60.00$ & $50.00$ & $420.00$ \\
Nielsen et al. \cite{nielsen_heat_1996} & $1.2306$ & $45.00$ & $36.00$ & $434.00$\\
Khorasany \& Fesanghary \cite{khorasany_novel_2009} & $1.1895$ & $28.10$ & $18.10$ & $451.91$ \\
this work - base & $1.1792$ & $25.52$ & $15.52$ & $454.48$ \\
this work - var UC & $1.1767$ & $25.59$ & $15.59$ & $454.41$ \\ \bottomrule
\end{tabular}
\end{table}

Implementing the CUS, we can find a solution with \mbox{$\SI{1.08}{\percent}$} lower TAC than the best literature value provided by Khorasany \& Fesanghary \cite{khorasany_novel_2009}. Compared to the base case, \mbox{$\SI{0.29}{\percent}$} can be saved. The stream matches of the two resulting HENs are both identical, see \ref{sec:appendixCaseStudyPlotsCS1}. The outlet temperature of the CUS of \mbox{$\SI{31.5}{\celsius}$} is very close to the lower valid range of \mbox{$\SI{31}{\celsius}$}. In contrast to the CU of the base case, the CUS temperature difference between inlet and outlet is reduced from \mbox{$\SI{15}{\kelvin}$} to \mbox{$\SI{1.5}{\kelvin}$}. The reduced temperature difference results in a large LMTD, which in turn results in a smaller and less expensive HEX area.

\subsection{Case study 2}
In the second case study, a frequently discussed aromatics plant in the literature is considered. The stream data provided by Linnhoff \& Ahmad \cite{linnhoff_cost_1990} is given in Table \ref{tab:UseCase2}. In this case study, thermal oil is used as hot utility. Since only the sensible heat of the oil is used, the outlet temperature can be adjusted and the hot utility is implemented as HUS. The cold utility uses the sensible heat of cooling water and is therefore implemented as CUS. In contrast to CS1, two streams with variable outlet temperatures and flow capacities are implemented. 

\begin{table}[H]
    \centering
    \caption{Stream data for case study 2: Linnhoff \& Ahmad \cite{linnhoff_cost_1990}.}
    \label{tab:UseCase2}
    \begin{threeparttable}
    \begin{tabular}{crrrr}
        \toprule
        Stream & $T^{\mathrm{in}}$ ($\SI{}{\celsius}$) & $T^{\mathrm{out}}$ ($\SI{}{\celsius}$) & $F$ ($\SI{}{\kilo\watt\per\kelvin}$) & $h$ ($\SI{}{\kilo\watt\per\meter\squared\per\kelvin}$)\\
        \midrule 
        H1 & $327$ & $40$ & $100$ & $0.50$ \\
        H2 & $220$ & $160$ & $160$ & $0.40$ \\
        H3 & $220$ & $60$ & $60$ & $0.14$ \\
        H4 & $160$ & $45$ & $400$ & $0.30$ \\
        C1 & $100$ & $300$ & $100$ & $0.35$ \\
        C2 & $35$ & $164$ & $70$ & $0.70$ \\
        C3 & $85$ & $138$ & $350$ & $0.50$ \\
        C4 & $60$ & $170$ & $60$ & $0.14$ \\
        C5 & $140$ & $300$ & $200$ & $0.60$ \\
        UH & $330$ & $250$ & - & $0.50$ \\
        UH$^\text{v}$ & $330$ & $[329, 250]$ & $(0, 25000]$ & $0.50$ \\
        UC & $15$ & $30$ & - & $0.50$ \\
        UC$^\text{v}$ & $15$ & $[16, 30]$ & $(0, 35000]$ & $0.50$ \\
        \bottomrule
    \end{tabular}
    \begin{tablenotes}[param, flushleft]
    \footnotesize
        \item HEX costs: \mbox{$c_{\mathrm{f}} = \SI{2000}{\USD\per\year}$}, \mbox{$c_{\mathrm{v}} = \SI{70}{\USD\per\betaCo\per\year}$}, \mbox{$\beta = 1$}
        \item Utility costs: \mbox{$c_{\mathrm{hu}} =  \SI{60}{\USD\per\kilo\watt\per\year}$}, \mbox{$c_{\mathrm{cu}} =  \SI{6}{\USD\per\kilo\watt\per\year}$}
        \item Min. approach temperature: \mbox{$\Delta T_{\mathrm{min}} = \SI{1}{\celsius}$}
    \end{tablenotes}
    \end{threeparttable}
\end{table}

As can be seen from Table \ref{tab:performanceDataCS2}, the number of binary variables for the var UC \& UH case has more than tripled compared to the base case. Due to the greater complexity, the computation time also increases from \mbox{$\SI{133.95}{\second}$} to \mbox{$\SI{590.12}{\second}$}.

\begin{table}[htp]
\centering
\caption{Problem size, computation time and relative gap for CS2.}
\label{tab:performanceDataCS2}
\begin{tabular}{lrrrrr}
\toprule
\multicolumn{1}{c}{case} & \multicolumn{1}{c}{variables / -} & \multicolumn{1}{c}{binaries / -}  & \multicolumn{1}{c}{time / s} & \multicolumn{1}{c}{rel. gap / \%} \\
\midrule
 base & $1009$ & $387$ & $133.95$ & $0.0096$ \\
 var UC \& UH & $3014$ & $1201$ & $590.12$ & $0.0000$ \\
 \bottomrule
\end{tabular}
\end{table}

Table \ref{tab:results-UseCase2} lists a comparison of the calculated values with data from the literature. Linnhoff \& Ahmad \cite{linnhoff_cost_1990} used the pinch design method and the driving force plot to obtain a HEN with \mbox{$2.9300 \cdot 10^6 \, \SI{}{\USD\per\year}$} of TAC. Through evolution and continuous optimization of the exchanger duties, they were able to obtain TAC of \mbox{$2.8900 \cdot 10^6 \, \SI{}{\USD\per\year}$}. In both calculations, the minimum temperature difference at the hot utility of \mbox{$\SI{26}{\celsius}$} was violated, causing the outlet temperature of the thermal oil to be higher than \mbox{$\SI{250}{\celsius}$}. Fieg et al. \cite{fieg_monogenetic_2009} corrected this by calculating the hot utility costs proportional to the thermal oil mass flow rate as a function of the outlet temperature. The TAC of \mbox{$2.9300 \cdot 10^6 \, \SI{}{\USD\per\year}$} were corrected to \mbox{$2.9920 \cdot 10^6 \, \SI{}{\USD\per\year}$}, and of \mbox{$2.8900 \cdot 10^6 \, \SI{}{\USD\per\year}$} to \mbox{$3.0250 \cdot 10^6 \, \SI{}{\USD\per\year}$}, respectively. Fieg et al. \cite{fieg_monogenetic_2009} also calculated an even cheaper solution with \mbox{$2.9223 \cdot 10^6 \, \SI{}{\USD\per\year}$} of TAC using a hybrid genetic algorithm (GA). Lewin \cite{lewin_generalized_1998} used a GA to find several solutions using different parameters of the algorithm. Despite an assumed minimum temperature difference of \mbox{$\SI{10}{\celsius}$}, the outlet temperature of the hot utility was raised without considering mass flow dependent costs. The best solution has TAC of \mbox{$2.9360 \cdot 10^6 \, \SI{}{\USD\per\year}$}. Zhu et al. \cite{zhu_method_1995} found the most expensive solution so far with TAC of \mbox{$2.9700 \cdot 10^6 \, \SI{}{\USD\per\year}$} by a two-step procedure using heuristics and nonlinear optimization. The optimizations performed in the literature enable a HEN with three stages and stream splits. In this paper, we assume on the one hand that the hot utility costs are proportional to the heat flow and not to the mass flow. On the other hand, we use a lower bound for the minimum temperature difference of \mbox{$\SI{1}{\celsius}$}. For the base case, TAC of \mbox{$2.9114 \cdot 10^6 \, \SI{}{\USD\per\year}$} could be calculated using only two stages and stream splits. Compared to the best literature value from Linnhoff \& Ahmad \cite{linnhoff_cost_1990}, this solution is \mbox{$\SI{0.74}{\percent}$} more expensive. Compared to the worst literature value from Zhu et al. \cite{zhu_method_1995}, however, it is \mbox{$\SI{2.01}{\percent}$} cheaper. 

\begin{table}[htp]
\centering
\caption{Results for CS2: Comparison of costs and heat loads at the utilities and streams.}
\label{tab:results-UseCase2}
\begin{threeparttable}
\begin{tabular}{lcccc}
\toprule
\multirow{2}{*}{reference} &  & \multicolumn{3}{c}{heat load / $\SI{}{\mega\watt}$} \\
 & \multicolumn{1}{c}{TAC / $10^6 \, \SI{}{\USD\per\year}$} & \multicolumn{1}{c}{CU} & \multicolumn{1}{c}{HU} & \multicolumn{1}{c}{stream} \\ \midrule
Linnhoff \& Ahmad \cite{linnhoff_cost_1990}\tnote{a} & $2.9300$ & $32.76$ & $25.04$ & $61.14$ \\
Linnhoff \& Ahmad \cite{linnhoff_cost_1990}\tnote{b} & $2.8900$ & $33.03$ & $25.31$ & $60.87$ \\
Zhu et al. \cite{zhu_method_1995} & $2.9700$ & $33.94$ & $26.22$ & \tnote{/*} \\
Lewin \cite{lewin_generalized_1998} & $2.9360$ & $32.81$ & $25.09$ & $61.09$ \\
Fieg et al. \cite{fieg_monogenetic_2009} & $2.9223$ & $31.34$ & $23.62$ & $62.57$ \\
this work - base & $2.9114$ & $31.42$ & $23.70$ & $62.48$ \\
this work - var UC \& UH & $2.8526$ & $31.43$ & $23.71$ & $62.47$ \\ \bottomrule
\end{tabular}
\begin{tablenotes}
    \item[a] Figure 19 (a) of \cite{linnhoff_cost_1990}.
    \item[b] Figure 19 (b) of \cite{linnhoff_cost_1990}.
    \item[/*] Heat load is not reported by the authors.
  \end{tablenotes}
\end{threeparttable}
\end{table}

The TAC can be reduced to \mbox{$2.8526 \cdot 10^6 \, \SI{}{\USD\per\year}$} by implementing a CUS and a HUS. Compared to the best literature value provided by Linnhoff \& Ahmad \cite{linnhoff_cost_1990}, the solution of the var UC \& UH case is \mbox{$\SI{1.30}{\percent}$} cheaper. The heat loads differ only slightly in both cases. The stream plots are shown in \ref{sec:appendixCaseStudyPlotsCS2}. The stream matches of the hot and cold streams are identical in both cases. The hot utility is used on the cold streams C1, C2, and C5 in both cases. The outlet temperature of the hot utility increases from \mbox{$\SI{250.0}{\celsius}$} to \mbox{$\SI{328.1}{\celsius}$}. The cold utility is used in the base case and the var UC \& UH case on streams H1, H3 and H4. In the var UC \& UH case, the cold utility is used a second time on stream H4. In this case, the CUS exchanges heat at two stages with the hot streams. The outlet temperature of the cooling water decreases from \mbox{$\SI{30.0}{\celsius}$} to \mbox{$\SI{16.8}{\celsius}$}. Again, in this case study, the cost savings are due to reduced HEX areas due to a higher thermal oil outlet temperature and a lower cooling water outlet temperature.

\subsection{Case study 3}
The third case study is the Bandar Imam aromatic plant on the northwestern coast of the Persian Gulf. The real-world problem includes six hot and ten cold streams. The stream data was provided by Khorasany \& Fesanghary \cite{khorasany_novel_2009} and is listed in \mbox{Table \ref{tab:UseCase3}}. Khorasany \& Fesanghary \cite{khorasany_novel_2009} probably made a conversion error from $\SI{}{\celsius}$ to $\SI{}{\kelvin}$ for the cold utility  \cite{pavao_new_2018, feyli_reliable_2022, nair_unified_2019, kayange_non-structural_2020}. Since cooling water is used as cold utility, it is implemented as CUS. In this case, two hot utilties are available: flue gas (UH$_1$) and steam (UH$_2$). The flue gas utility is implemented as HUS. The steam utility is implemented as a conventional utility without variable temperatures and flow capacity.

\begin{table}[htp]
    \centering
    \caption{Stream data for case study 3: Khorasany \& Fesanghary \cite{khorasany_novel_2009}.}
    \label{tab:UseCase3}
    \begin{threeparttable}
    \begin{tabular}{crrrr}
        \toprule
        Stream & $T^{\mathrm{in}}$ ($\SI{}{\celsius}$) & $T^{\mathrm{out}}$ ($\SI{}{\celsius}$) & $F$ ($\SI{}{\kilo\watt\per\kelvin}$) & $h$ ($\SI{}{\kilo\watt\per\meter\squared\per\kelvin}$)\\
        \midrule 
        H1 & $385.0$ & $159.0$ & $131.51$ & $1.238$ \\
        H2 & $516.0$ & $43.0$ & $1198.96$ & $0.546$ \\
        H3 & $132.0$ & $82.0$ & $378.96$ & $0.771$ \\
        H4 & $91.0$ & $60.0$ & $589.55$ & $0.859$ \\
        H5 & $217.0$ & $43.0$ & $186.22$ & $1.000$ \\
        H6 & $649.0$ & $43.0$ & $116.00$ & $1.000$ \\
        C1 & $30.0$ & $385.0$ & $119.10$ & $1.850$ \\
        C2 & $99.0$ & $471.0$ & $191.05$ & $1.129$ \\
        C3 & $437.0$ & $521.0$ & $377.91$ & $0.815$ \\
        C4 & $78.0$ & $418.6$ & $160.43$ & $1.000$ \\
        C5 & $217.0$ & $234.0$ & $1297.70$ & $0.443$ \\
        C6 & $256.0$ & $266.0$ & $2753.00$ & $2.085$ \\
        C7 & $49.0$ & $149.0$ & $197.39$ & $1.000$ \\
        C8 & $59.0$ & $163.4$ & $123.56$ & $1.063$ \\
        C9 & $163.0$ & $649.0$ & $95.98$ & $1.810$ \\
        C10 & $219.0$ & $221.3$ & $1997.50$ & $1.377$ \\
        UH$_1$ & $1800.0$ & $800.0$ & - & $1.200$ \\
        UH$^\text{v}_1$ & $1800.0$ & $[800.0, 1700.0]$ & $(0, 150.00]$ & $1.200$ \\
        UH$_2$ & $509.0$ & $509.0$ & - & $1.000$ \\
        UC & $38.0$ & $82.0$ & - & $1.000$ \\
        UC$^\text{v}$ & $38.0$ & $[39.0, 82.0]$ & $[10000, 500000]$ & $1.000$ \\
        \bottomrule
    \end{tabular}
    \begin{tablenotes}[param, flushleft]
    \footnotesize
        \item HEX costs: \mbox{$c_{\mathrm{f}} = \SI{26600}{\USD\per\year}$}, \mbox{$c_{\mathrm{v}} = \SI{4147.5}{\USD\per\betaCo\per\year}$}, \mbox{$\beta = 0.6$}
        \item Hot utility costs: \mbox{$c_{\mathrm{hu,1}} =  \SI{35}{\USD\per\kilo\watt\per\year}$}, \mbox{$c_{\mathrm{hu,2}} =  \SI{27}{\USD\per\kilo\watt\per\year}$}
        \item Cold utility costs: \mbox{$c_{\mathrm{cu}} =  \SI{2.1}{\USD\per\kilo\watt\per\year}$}
        \item Min. approach temperature: \mbox{$\Delta T_{\mathrm{min}} = \SI{1}{\celsius}$}
    \end{tablenotes}
    \end{threeparttable}
\end{table}

Table \ref{tab:performanceDataCS3} lists the number of variables, computation time, and relative gap. Both cases were optimally solved within the defined MIP gap. Due to the high complexity of the var UC \& UH case, the computation time is the highest with \mbox{$\SI{11335.64}{\second}$} (\mbox{$\SI{3.15}{\hour}$}).

\begin{table}[htp]
\centering
\caption{Problem size, computation time and relative gap for CS3.}
\label{tab:performanceDataCS3}
\begin{tabular}{lrrrrr}
\toprule
\multicolumn{1}{c}{case} & \multicolumn{1}{c}{variables / -} & \multicolumn{1}{c}{binaries / -}  & \multicolumn{1}{c}{time / s} & \multicolumn{1}{c}{rel. gap / \%} \\
\midrule
 base & $2893$ & $1128$ & $637.07$ & $0.0094$ \\
 var UC \& UH & $7302$ & $2914$ & $\SI{11335.64}{}$ & $0.0100$ \\
 \bottomrule
\end{tabular}
\end{table}

Table \ref{tab:results-UseCase3} compares the TAC and heat loads of the existing power plant with those from literature and those calculated in this paper. Khorasany \& Fesanghary \cite{khorasany_novel_2009} used the same two-level approach as for CS1 to optimize the HEN of the aromatic plant. The optimized HEN with TAC of \mbox{$7.4357 \cdot 10^6 \, \SI{}{\USD\per\year}$} results in a cost savings of \mbox{$\SI{16.00}{\percent}$} compared to the existing plant. Feyli et al. \cite{feyli_reliable_2022} found TAC of \mbox{$7.1285 \cdot 10^6 \, \SI{}{\USD\per\year}$} using GA and a quasi-linear programming method. The HEN has eight stages and no stream splits. Aguitoni et al. \cite{aguitoni_heat_2018} found a slight improvement in TAC. Here the discrete variables were optimized with a GA and the heat loads on the streams, and the stream split fractions with the help of differential evolution. The resulting HEN has three stages with stream splits. Pavão et al. followed a two-level approach in \cite{pavao_enhanced_2018} and \cite{pavao_new_2018}, handling the binary variables with simulated annealing. The continuous variables are optimized using a rocket fireworks approach. In \cite{pavao_enhanced_2018}, the super-structure formulation was adapted, allowing intermediate placement of utilities. In \cite{pavao_new_2018}, sub-stages, sub-splits and cross-flows are additionally enabled. Nair \& Karimi \cite{nair_unified_2019} used a stageless superstructure formulation with stream splits to optimize the TAC. Here, MILP relaxations are solved to find a lower bound on the TAC. Constraining the configurations and solving the nonlinear problem yields to an upper bound on the TAC. The lowest TAC to date of \mbox{$6.6647 \cdot 10^6 \, \SI{}{\USD\per\year}$} could be determined by Liu et al. \cite{liu_extended_2022} using a GA approach. The HEN is composed of three stages and several stream splits. In this paper, we obtained TAC of \mbox{$6.7451 \cdot 10^6 \, \SI{}{\USD\per\year}$} for the base case with two stages and stream splits. The stream plots are shown in \ref{sec:appendixCaseStudyPlotsCS3}. Like most solutions from the literature, the second hot utility (steam) is not used. Compared to the best literature value of Liu et al. \cite{liu_extended_2022}, this solution is \mbox{$\SI{1.19}{\percent}$} more expensive. However, our solution is \mbox{$\SI{31.30}{\percent}$} cheaper than the existing plant.

\begin{table}[htp]
\centering
\caption{Results for CS3: Comparison of costs and heat loads at the utilities and streams.}
\label{tab:results-UseCase3}
\begin{threeparttable}
\begin{tabular}{lcccc}
\toprule
\multirow{2}{*}{reference} &  & \multicolumn{3}{c}{heat load / $\SI{}{\mega\watt}$} \\
 & \multicolumn{1}{c}{TAC / $10^6 \, \SI{}{\USD\per\year}$} & \multicolumn{1}{c}{CU} & \multicolumn{1}{c}{HU} & \multicolumn{1}{c}{stream} \\ \midrule
 existing plant \cite{khorasany_novel_2009} & $8.8564$ & $524.72$ & $122.16$ & \tnote{/*} \\
Khorasany \& Fesanghary \cite{khorasany_novel_2009} & $7.4357$ & $469.62$ & $66.07$ & $267.09$ \\
Feyli et al. \cite{feyli_reliable_2022} & $7.1285$ & $437.77$ & $34.21$ & $298.96$ \\
Aguitoni et al. \cite{aguitoni_heat_2018} & $7.1028$ & $437.44$ & $33.87$ & $299.29$ \\
Pavão et al. \cite{pavao_enhanced_2018} & $6.8013$ & $414.03$ & $10.47$ & $322.69$ \\
Pavão et al. \cite{pavao_new_2018} & $6.7126$ & $413.07$ & $9.50$ & $323.65$ \\
Nair \& Karimi \cite{nair_unified_2019} & $6.6956$ & $412.25$ & $8.69$ & $324.47$ \\
Liu et al. \cite{liu_extended_2022} & $6.6647$ & $413.11$ & $9.55$ & $323.61$ \\
this work - base & $6.7451$ & $414.72$ & $11.13$ & $322.04$ \\
this work - var UC \& UH & $6.3297$ & $416.57$ & $12.98$ & $320.18$ \\ \bottomrule
\end{tabular}
\begin{tablenotes}
    \item[/*] Heat load is not reported by the authors.
  \end{tablenotes}
\end{threeparttable}
\end{table}

In the var UC \& UH case, the TAC can be reduced to \mbox{$6.3297 \cdot 10^6 \, \SI{}{\USD\per\year}$} through the implementation of a HUS and a CUS. Thus, this solution is \mbox{$\SI{5.37}{\percent}$} cheaper than the best solution provided by Liu et al. \cite{liu_extended_2022}. Compared to the existing plant; the TAC can be reduced by \mbox{$\SI{40.49}{\percent}$}. In the base and var UC \& UH case, the second utility (steam) is not used. The HEN configuration of the two cases differs only in the \nth{1} and \nth{5} hot stream, see \ref{sec:appendixCaseStudyPlotsCS3}. The outlet temperature of the first utility (flue gas) is increased from \mbox{$\SI{800.0}{\celsius}$} to \mbox{$\SI{1674.4}{\celsius}$}. The outlet temperature of the cold utility is lowered from \mbox{$\SI{82.0}{\celsius}$} to \mbox{$\SI{39.6}{\celsius}$}. Therefore, the LMTD at the HEX will be increased, resulting in smaller and more cost-efficient HEX areas.
\section{Conclusion}
\label{sec:conclusion}

We have adapted the superstructure formulation from Yee \& Grossmann \cite{yee_simultaneous_1990} to implement utilities as streams with variable temperatures and flow capacity. All non-linear correlations, such as the LMTD, reduced HEX areas and energy balances of the streams are piecewise-linear approximated with simplex and hyperplane models. The translation to MILP is achieved highly efficiently using logarithmic coding. A possible application of this method is the optimal design of utilities. Especially for utilities, which use only sensible heat from, for example, thermal oil, water or flue gas, it is possible to adjust the temperatures within defined limits. The implementation of hot and cold utilities as streams also allows a cost-efficient multi-stage heat exchange with stream splits.

The methods presented in this paper were applied to three representative case studies with 4, 9 and 16 streams, respectively. To verify the developed framework, a base case HEN without utility streams was calculated and compared with results from the literature for all case studies. For the first and second case study, the best available solutions from the literature were improved in terms of TAC by \mbox{$\SI{0.87}{\percent}$} and \mbox{$\SI{0.74}{\percent}$}, respectively. For the third case study, TAC which is only \mbox{$\SI{1.19}{\percent}$} higher than the best literature value was calculated. Since the base case results are similar to the literature values, it can be assumed that the HENS formulation with piecewise-linear approximations of the non-linear terms provides reliable results. For CS1, cost savings of \mbox{$\SI{0.29}{\percent}$} can be achieved by implementing the cold utility as a stream. For CS2 and CS3, both the hot and cold utilities were implemented as streams. Cost reductions of \mbox{$\SI{1.30}{\percent}$} and \mbox{$\SI{5.37}{\percent}$}, respectively, were achieved. The results show that in a HEN with utilities implemented as a stream, the outlet temperature of cold utilities should always be as low as possible and that of hot utilities should always be as high as possible. The resulting small temperature differences between inlet and outlet temperatures at the utilities result in the largest possible LMTD at the heat exchangers. This reduces the HEX area for the same amount of heat to be transferred, which in turn reduces the TAC. However, it should be mentioned that a small temperature difference between utility stream supply and return leads to high flow capacities. Compared to the utilities of the base case, the mass flows must be increased to achieve a smaller temperature difference. The assumption of mass flow independent utility costs is only valid if the required pumping power is negligible compared to the costs of heat input and output. 

In summary, our research indicates the need to implement utilities as streams with variable outlet temperatures and flow capacities. By setting a variable outlet temperature, where it does not have to be necessarily defined a priori, a previously untapped potential for cost and energy savings is created. Additional savings potential arises from the possibility of multi-stage heat exchange and stream splits of the utility streams. In this way, energy-saving methods can be more cost-effective, encouraging their widespread use. This research also provides a foundation for the MILP implementation of streams with variable temperatures and/or flow capacities. Thus, in future work, the operational behavior of plants that effect stream parameters can be integrated into HENS. This will open up new possibilities for energy- and cost-optimized holistic process optimization.

\bibliographystyle{elsarticle-num-names} 
\bibliography{_main_paper}
%\bibliography{bibliography}

\section*{Statements and Declarations}
\subsection*{Funding}
The work this paper is based on was funded by the Austrian research promotion agency (FFG) under grant number 884340. 

The authors acknowledge TU Wien Bibliothek for financial support for editing/proofreading.

\subsection*{Competing Interests}
The authors have no relevant financial or non-financial interests to disclose.

\subsection*{Authors Contributions}
The method presented in this paper was developed by David Huber. Testing and evaluation was conducted by David Huber and Felix Birkelbach. The conceptualization of the paper was the work of all authors. The first draft was written by David Huber. All authors contributed to the revision of the initial draft. René Hofmann was responsible for funding and supervision. All authors read and approved the final manuscript.

\appendix
\section{Convergence Behavior}
\label{sec:appendixConvergence}
\begin{figure}[H]
    \centering    \includegraphics[width=1\textwidth]{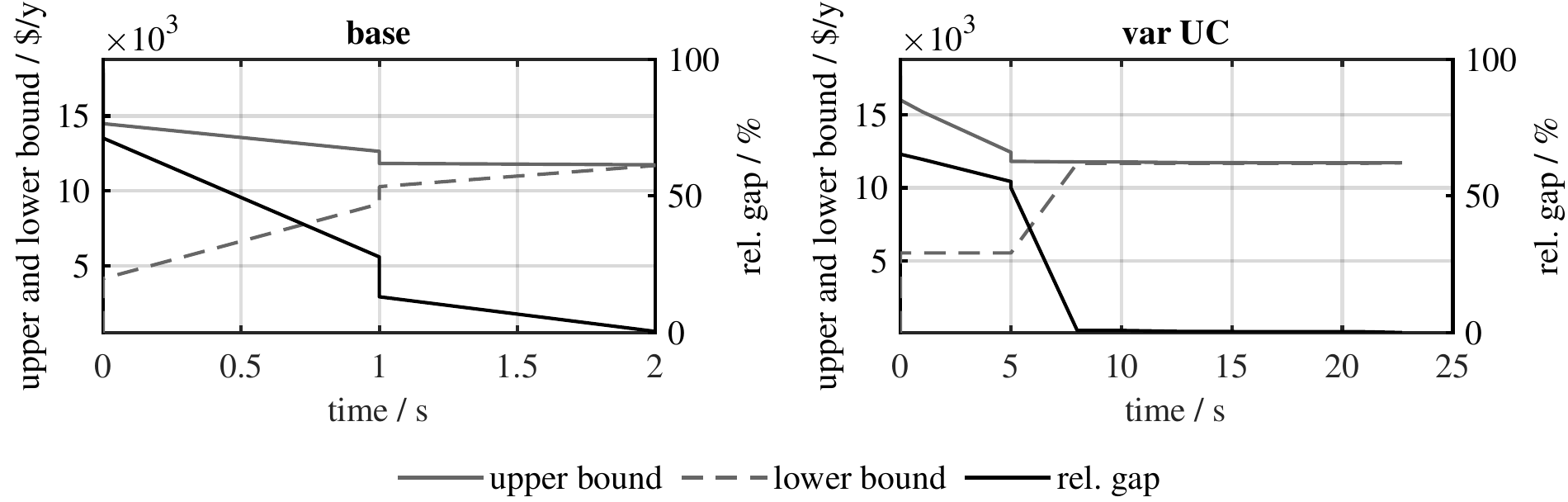}
    \caption{Convergence behavior of CS1 for the base and the var UC case.}
    \label{fig:convergenceCS1}
\end{figure}

\begin{figure}[H]
    \centering    
    \includegraphics[width=1\textwidth]{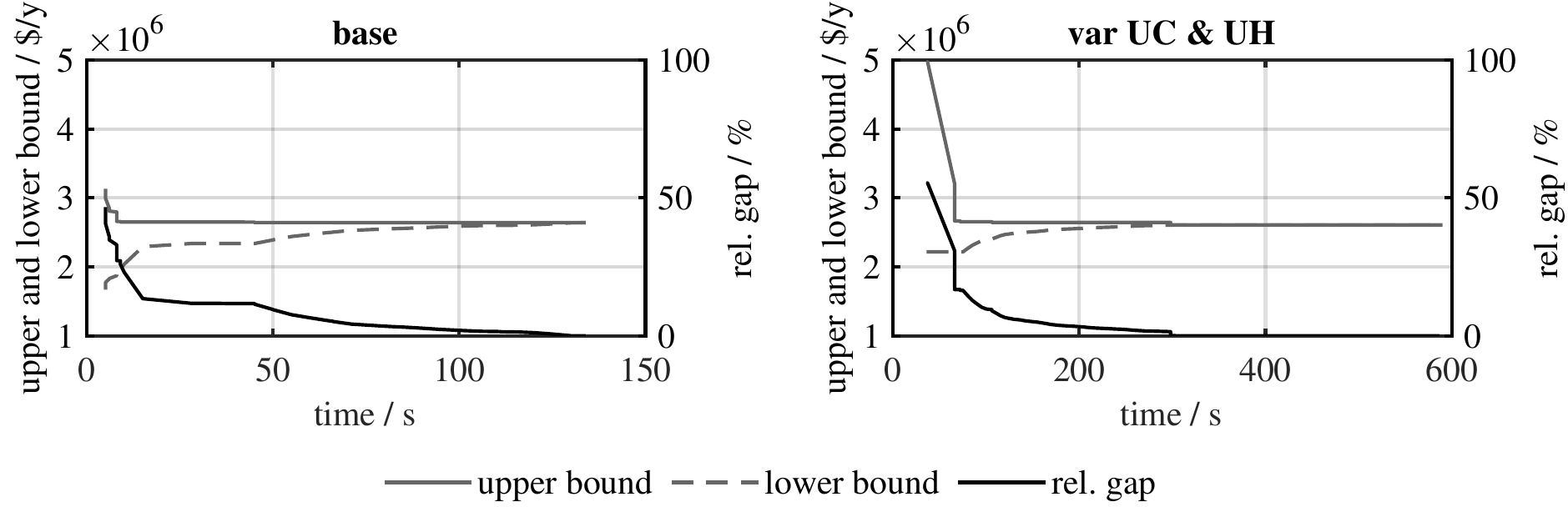}
    \caption{Convergence behavior of CS2 for the base and the var UC \& UH case.}
    \label{fig:convergenceCS2}
\end{figure}

\begin{figure}[H]
    \centering    
    \includegraphics[width=1\textwidth]{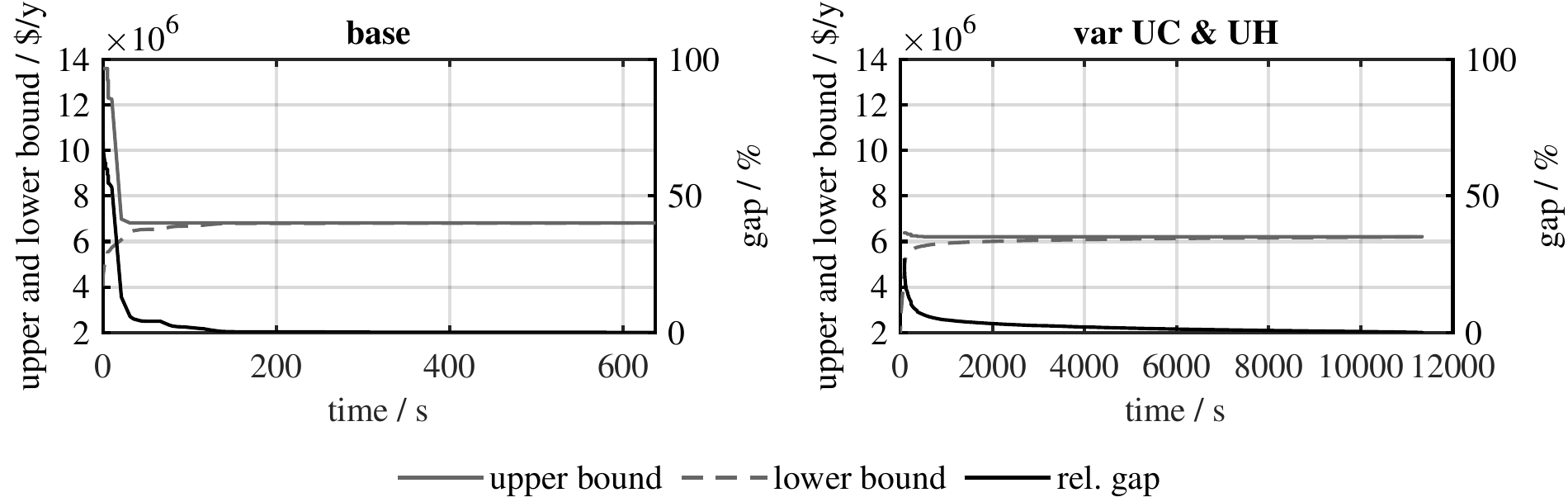}
    \caption{Convergence behavior of CS3 for the base and the var UC \& UH case.}
    \label{fig:convergenceCS3}
\end{figure}

The convergence behavior of the three case-studies is based on the log-file of the solver and is shown in Figures \ref{fig:convergenceCS1} to \ref{fig:convergenceCS3}. The script to process the unformatted data into processable vectors was coded by ChatGPT \cite{assistant_openai_2022}.

\subsection{Convergence Behavior CS1}

\subsection{Convergence Behavior CS2}

\subsection{Convergence Behavior CS3}

\section{Stream Plots}
\label{sec:appendixCaseStudyPlots}
\begin{figure}[H]
    \centering
    \includegraphics[width=1\textwidth]{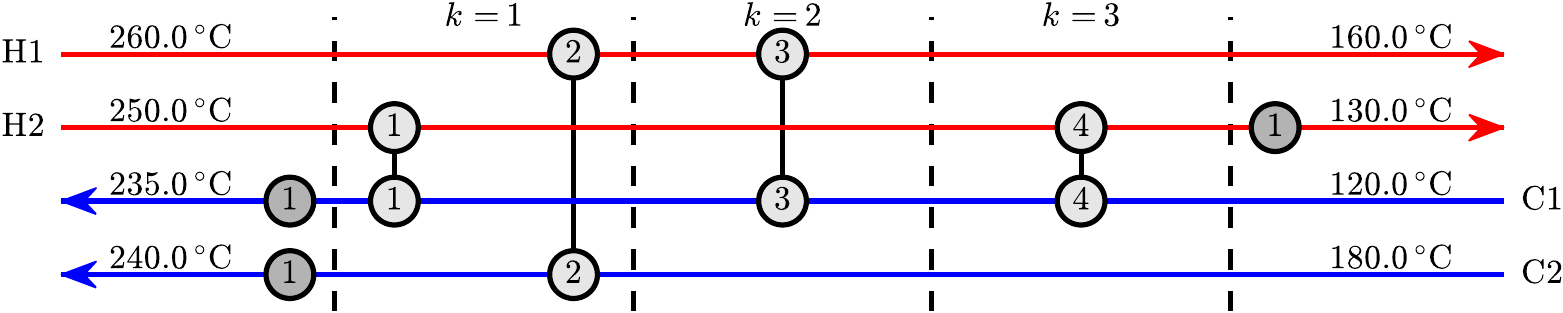}
    \caption{Optimized HEN configuration of CS1: base case with TAC of \mbox{$1.1792 \cdot 10^4 \, \SI{}{\USD\per\year}$}.}
    \label{fig:case1Base}
\end{figure}

\begin{figure}[H]
    \centering
    \includegraphics[width=1\textwidth]{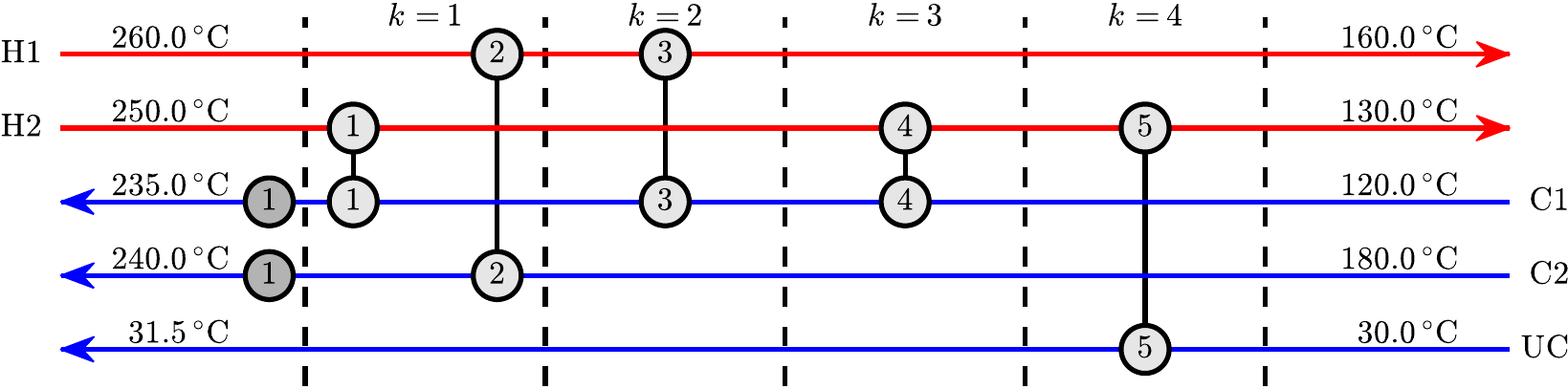}
    \caption{Optimized HEN configuration of CS1: var UC case with TAC of \mbox{$1.1767 \cdot 10^4 \, \SI{}{\USD\per\year}$}.}
    \label{fig:case1VarUC_F_Tout}
\end{figure}

\begin{figure}[H]
    \centering
    \includegraphics[width=1\textwidth]{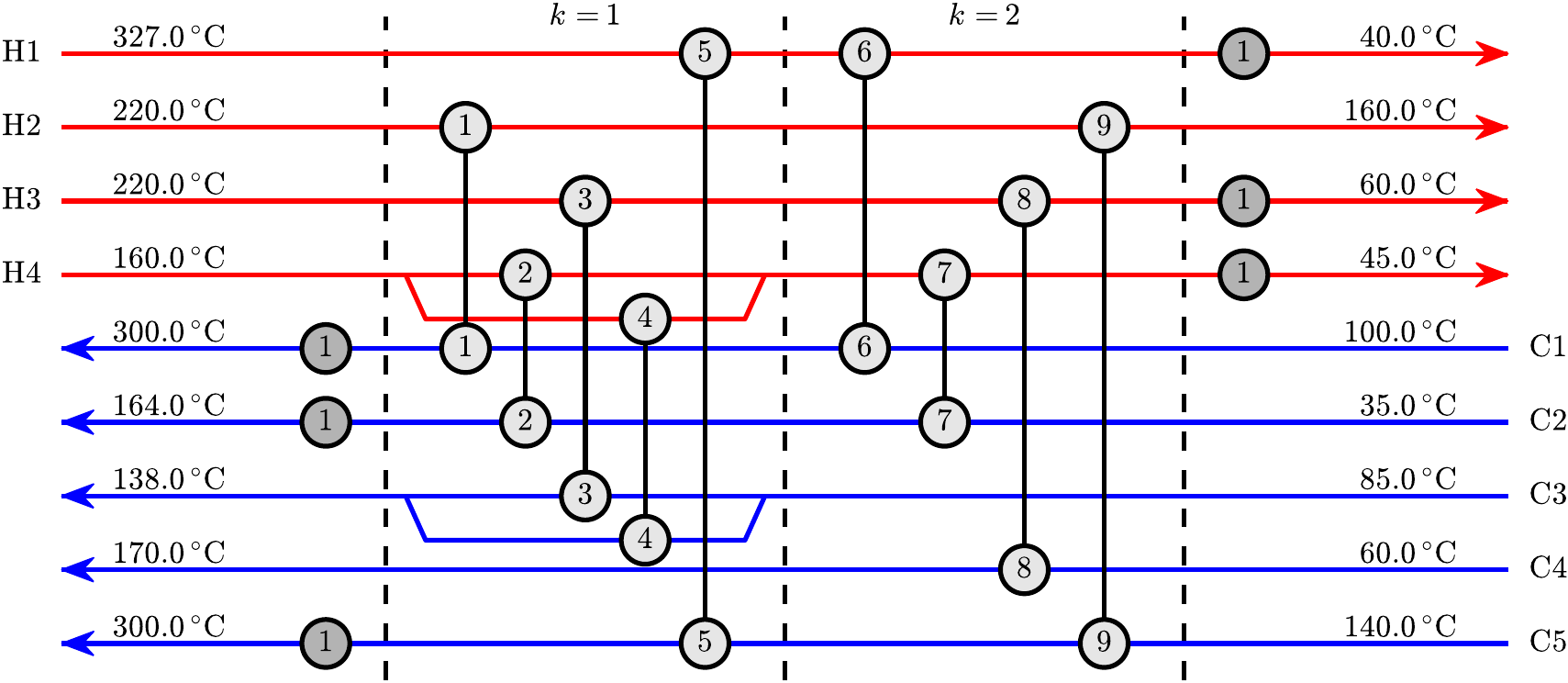}
    \caption{Optimized HEN configuration of CS2: base case with TAC of \mbox{$2.9114 \cdot 10^6 \, \SI{}{\USD\per\year}$}.}
    \label{fig:case2Base}
\end{figure}

\begin{figure}[H]
    \centering
    \includegraphics[width=1\textwidth]{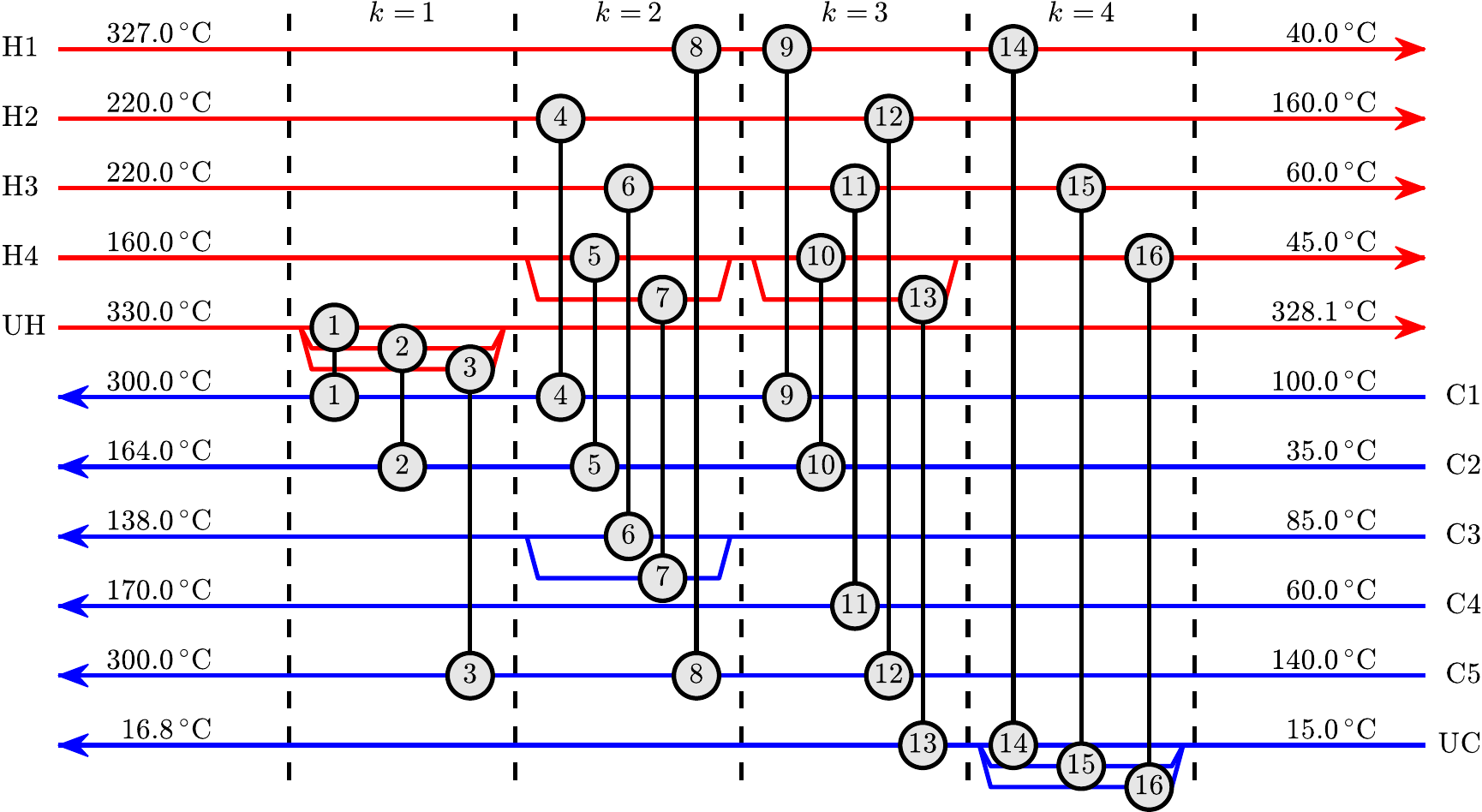}
    \caption{Optimized HEN configuration of CS2: var UC \& UH case with TAC of \mbox{$2.8526 \cdot 10^6 \, \SI{}{\USD\per\year}$}.}
    \label{fig:case2varUCUH}
\end{figure}

\begin{figure}[H]
    \centering
    \includegraphics[width=1\textwidth]{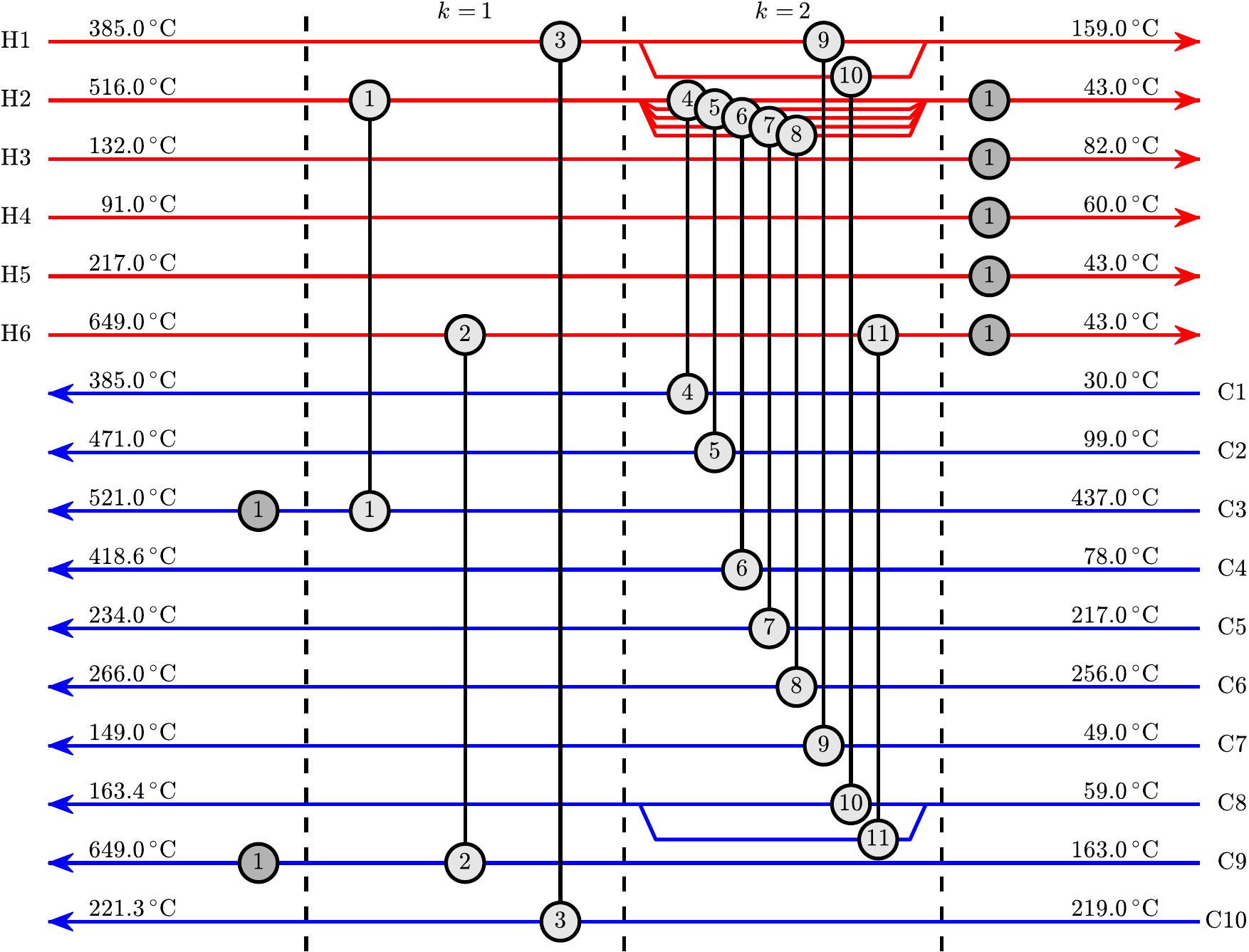}
    \caption{Optimized HEN configuration of CS3: base case with \mbox{$6.7451 \cdot 10^6 \, \SI{}{\USD\per\year}$}. Note that in this case only the flue gas is used as a hot utility.}
    \label{fig:case3Base}
\end{figure}

\begin{figure}[H]
    \centering
    \includegraphics[width=1\textwidth]{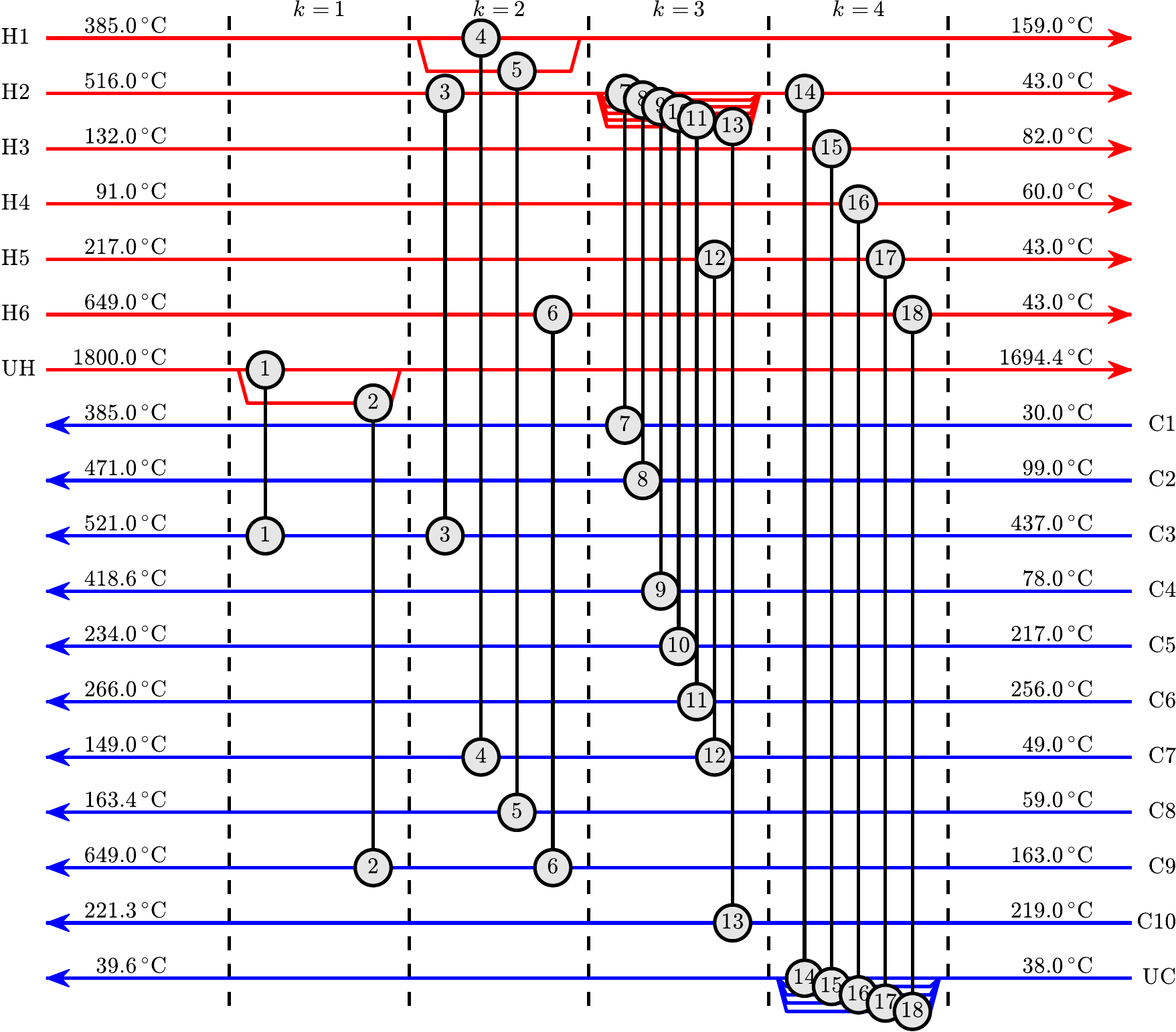}
    \caption{Optimized HEN configuration of CS3: var UC \& UH case with TAC of \mbox{$6.3297 \cdot 10^6 \, \SI{}{\USD\per\year}$}. Note that in this case only the flue gas is used as a hot utility.}
    \label{fig:case3varUCUH}
\end{figure}

Figures \ref{fig:case1Base} to \ref{fig:case3varUCUH} show the optimized stream plots of the three case studies from Section \ref{sec:usecase}. Hot streams are shown in red and cold streams in blue. The light gray circles inside the stages $k$ represent heat exchangers. The dark gray circles without connecting lines are hot and cold utilities, respectively.

\subsection{Stream Plots CS1}
\label{sec:appendixCaseStudyPlotsCS1}

\subsection{Stream Plots CS2}
\label{sec:appendixCaseStudyPlotsCS2}

\subsection{Stream Plots CS3}
\label{sec:appendixCaseStudyPlotsCS3}

\end{document}